\documentclass[11pt,a4paper]{amsart}
\usepackage[a4paper,margin=0.8in,includefoot,footskip=0.45in]{geometry}
\usepackage{amssymb,graphicx}
\usepackage[T1]{fontenc}
\usepackage[utf8]{inputenc}
\usepackage{microtype}
\usepackage{enumitem}
\setlist[enumerate]{font=\normalfont}
\usepackage{needspace}
\usepackage[hidelinks]{hyperref}
\swapnumbers
\newtheoremstyle{boldheadplain}
  {1.2\baselineskip}{0.7\baselineskip}
  {\itshape}{}{\bfseries}{.}{.5em}{}
\newtheoremstyle{boldheadupright}
  {1.2\baselineskip}{0.7\baselineskip}
  {\normalfont}{}{\bfseries}{.}{.5em}{}

\theoremstyle{boldheadplain}
\newtheorem{theorem}{Theorem}[section]
\newtheorem{lemma}[theorem]{Lemma}
\newtheorem{proposition}[theorem]{Proposition}
\newtheorem*{theoremA}{Theorem A}
\newtheorem*{theoremB}{Theorem B}
\newtheorem{corollary}[theorem]{Corollary}
\theoremstyle{boldheadupright}\newtheorem{remark}[theorem]{Remark}\newtheorem{definition}[theorem]{Definition}
\AddToHook{env/definition/before}{\Needspace{4\baselineskip}}
\AddToHook{env/lemma/before}{\Needspace{4\baselineskip}}
\AddToHook{env/proposition/before}{\Needspace{4\baselineskip}}
\AddToHook{env/theorem/before}{\Needspace{4\baselineskip}}
\AddToHook{env/corollary/before}{\Needspace{4\baselineskip}}
\AddToHook{env/remark/before}{\Needspace{4\baselineskip}}
\newcommand{\Aut}{\operatorname{Aut}}
\newcommand{\PA}{\mathsf{PA}}
\newcommand{\ZF}{\mathsf{ZF}}
\newcommand{\Lset}{\mathcal L_{\mathrm{Set}}}
\newcommand{\Sk}{\operatorname{Sk}}
\newcommand{\G}{\mathcal G}
\newcommand{\Grp}{\mathsf{Grp}}
\newcommand{\Grpr}{\Grp_{\mathrm r}}
\allowdisplaybreaks[1]
\title[Automorphisms] {Automorphism groups of models of set theory}
\author[]{Ali Enayat}
\address{Ali Enayat, Department of Philosophy, Linguistics, and Theory of Science, University of Gothenburg, Sweden}
\email{ali.enayat@gu.se}
\date{\today}
\makeatletter
\def\section{\@startsection{section}{1}{\z@}%
  {1.5\baselineskip}{0.8\baselineskip}%
  {\normalfont\large\bfseries\scshape\centering}}
\def\subsection{\@startsection{subsection}{2}{\z@}%
  {1.2\baselineskip}{0.5\baselineskip}%
  {\normalfont\normalsize\bfseries\scshape}}
\makeatother

\begin{document}
\begin{abstract}
Assuming the consistency of \textup{\textsf{ZF}}, we study the groups that arise as full automorphism groups of models of set theory, both in the presence and in the absence of Foundation. Our principal results are as follows.
\medskip

\noindent\textbf{Theorem A.} \emph{Let $\Grp$ be the class of all groups and $\Grpr$ the class of right-orderable groups}. 

\medskip

\begin{enumerate}
\item \emph{ The groups realizable as automorphism groups of models of any of the theories $\mathsf{ZFC}$, $\ZF+\mathsf{LM}$, and $\ZF+\mathsf{V=OD}$ are precisely the members of $\Grpr$ (the corresponding class for \textup{\textsf{ZF}} is strictly larger). Here \textup{\textsf{LM}} is the Leibniz--Mycielski axiom.}

\medskip

\item \emph{If $\mathcal M$ is a model of $\ZF+\mathsf{V=OD}$ and $G\in\Grpr$, there is $\mathcal N\equiv\mathcal M$ with $\Aut(\mathcal N)\cong G$ such that the canonically interpreted arithmetic structure of $\mathcal M$ admits an elementary end embedding into that of $\mathcal N$, preserving all relations induced by parameter-free formulas of set theory.}
\end{enumerate}

\medskip

\noindent\textbf{Theorem B.} \emph{If $\kappa$ is strongly inaccessible, then for every group $G$ of cardinality less than $\kappa$, there is a model of cardinality $\kappa$ of $\mathsf{ZFC}$ with Foundation replaced by its negation, whose automorphism group is isomorphic to $G$ and whose externally well-founded part is isomorphic to $(V_\kappa,\in)$.}
\end{abstract}
\makeatletter
\let\draftoriginalsettitle\@settitle
\def\@settitle{%
  \begingroup
  \centering
  \setlength{\fboxsep}{8pt}%
  \fbox{\normalfont\scshape Preliminary Draft}\par
  \endgroup
  \vskip 1em\relax
  \draftoriginalsettitle}
\makeatother
\maketitle
\makeatletter
\begingroup
\small
\setlength{\parskip}{0pt}
\renewcommand{\l@section}{\@tocline{1}{0pt}{1pc}{}{\normalfont\scshape}}
\renewcommand{\tocsection}[3]{\indentlabel{\normalfont\@ifnotempty{#2}{\ignorespaces#1 #2.\quad}}#3}
\tableofcontents
\endgroup
\makeatother

\newpage

\section{Introduction}
\label{sec:intro}
Given a first-order theory $T$, what groups can arise as the full automorphism groups of models of $T$? We refer to this question as the \emph{realization problem} for $T$. In this paper we investigate the realization problem for several extensions of Zermelo--Fraenkel set theory, as well as for set theory without Foundation, drawing in part on the corresponding results for Peano arithmetic.\footnote{For the study of fixed point sets of individual automorphisms of models of set theory, see \cite{EM,EKM,Togha}.} Throughout the paper we assume that $\ZF$ is consistent.

For a theory $T$ we write $\G(T)$ for the class of all groups isomorphic to the automorphism group of some model of $T$, and refer to this class as the \emph{automorphism-group spectrum} of $T$. We write $\Grp$ for the class of all groups and $\Grpr$ for the class of all right-orderable groups, i.e., groups admitting a linear ordering invariant under right multiplication.

A weaker version of the realization problem asks which groups can be \emph{embedded} in the automorphism group of some model of a given theory. The classical theorem of Ehrenfeucht and Mostowski on indiscernibles provides a general result in this direction: if $T$ has an infinite model, then for every linear order $\mathcal I$ there is a model $\mathcal M$ of $T$ such that $\Aut(\mathcal I)$ embeds into $\Aut(\mathcal M)$, where $\Aut(\mathcal K)$ denotes the automorphism group of a structure $\mathcal K$. Consequently, every right-orderable group embeds into the automorphism group of some model of any first-order theory with an infinite model (see Theorem~\ref{orderability}). Arranging that the embedding be surjective, however, requires additional control over the structure of $\mathcal M$. 

For models of Peano arithmetic (\textup{\textsf{PA}}), as shown by a fundamental result of Gaifman\footnote{\cite[Theorems~4.10--4.11; see also pp.~223--224]{Gaifman}}, given any model $\mathcal M$ of $\PA$ and any linear order $\mathcal L$, there is an elementary end extension $\mathcal N$ of $\mathcal M$ with
\[
\Aut(\mathcal N/M)\cong\Aut(\mathcal L).
\]
Here $\Aut(\mathcal N/M)$ denotes the subgroup fixing every element of the universe $M$ of $\mathcal M$. Taking $\mathcal M$ to be the prime model of a completion of $\PA$ realizes $\Aut(\mathcal L)$ as the full automorphism group of a model of $\PA$. Schmerl \cite[Theorem~A]{Schmerl} subsequently established an optimal generalization of this result, allowing $\Aut(\mathcal L)$ to be replaced by any of its closed subgroups, with respect to the topology of pointwise convergence (see Subsection~\ref{subsec:auttop}). His work also yields a complete characterization of the groups realizable as full automorphism groups of models of $\PA$, stated below as Theorem~\ref{schmerl-spectrum}.

\begin{theorem}[{Schmerl \cite[Theorem~A]{Schmerl}}]\label{schmerl-spectrum}
A group is isomorphic to the automorphism group of a model of $\PA$ if and only if it is right-orderable. Thus
\[
\G(\PA)=\Grpr.
\]
\end{theorem}

Theorem~\ref{schmerl-spectrum} provides a natural point of departure for our investigation of automorphism groups of models of set theory. Indeed, Ackermann's coding allows us to relate the realization problem for arithmetic to that for finite set theory.\footnote{For a survey of the comparative model theory of arithmetic and set theory in the context of elementary extensions, see \cite{Bogota}.} In this coding \cite{Ack}, one associates to each model $\mathcal M$ of \textup{\textsf{PA}} a membership relation $E$ on its universe by stipulating that
\begin{equation}\label{eq:ack}
 a E b\quad\Longleftrightarrow\quad \mathcal M\models\text{``the $a$-th binary digit of $b$ is $1$''}.
\end{equation}
For the standard model of \textup{\textsf{PA}} the resulting structure is isomorphic to $(V_\omega,\in)$.

At the level of theories, however, this correspondence requires a qualification. Let $\ZF_{\mathrm{fin}}$ be obtained from \textup{\textsf{ZF}} by replacing Infinity by its negation, and let \textup{\textsf{TC}} assert that every set has a transitive closure. The role of the transitive closure axiom in the correspondence between arithmetic and finite set theory is brought out by the following theorem of Kaye and Wong. 

\begin{theorem}[{Kaye--Wong \cite{KW}}]
\label{kw-equivalence}
The theories $\PA$ and $\ZF_{\mathrm{fin}}+\mathsf{TC}$
are definitionally equivalent.\footnote{Definitional
equivalence (also known as synonymy) means that the theories have a common definitional extension. It is a strong form
of bi-interpretability; see \cite[Remark~2]{EL}.
Throughout this paper, bi-interpretability is understood to be parameter-free, since bi-interpretability with parameters need not induce an isomorphism between the automorphism groups of the structures involved. It is also worth noting that over elementary set
theory (Extensionality, Empty Set, Pairing, Union
and Replacement), $\mathsf{TC}$ is equivalent to
the scheme of $\in$-induction, by
\cite[Proposition~12]{KW}.
See also \cite[Remark~2.2(f)]{ESV}.}
\end{theorem}

Since parameter-free bi-interpretability induces isomorphisms between the automorphism groups of corresponding models, Theorems~\ref{schmerl-spectrum} and~\ref{kw-equivalence} yield the following corollary.

\begin{corollary}\label{finite-spectrum}
The automorphism-group spectra of $\PA$ and $\ZF_{\mathrm{fin}}+\mathsf{TC}$ coincide, and both consist precisely of the right-orderable groups. In other words,
\[
\G(\ZF_{\mathrm{fin}}+\mathsf{TC})=\G(\PA)=\Grpr.
\]
\end{corollary}

The qualification concerning \textup{\textsf{TC}} is essential; indeed, \textup{\textsf{PA}} and $\ZF_{\mathrm{fin}}$ are not even sententially equivalent\footnote{Sentential equivalence is a weaker form of bi-interpretability in which the compositions of the interpretations at work are required only to yield models elementarily equivalent to the original models; see \cite[Section 5]{ESV}.} \cite[Theorem 5.1]{ESV}.

The situation changes dramatically when the transitive closure axiom is dropped, as the following result, obtained in joint work of the author with Jim Schmerl and Albert Visser, shows. In this context, an $\omega$-model of $\ZF_{\mathrm{fin}}$ is a model all of whose sets are externally finite.

\begin{theorem}[{\cite[Corollary 4.4 and footnote 19]{ESV}}]\label{esv-spectrum}
Every group is isomorphic to the full automorphism group of an $\omega$-model of $\ZF_{\mathrm{fin}}$. In particular,
\[
\G(\ZF_{\mathrm{fin}})=\Grp.
\]
\end{theorem}

The contrast between Corollary~\ref{finite-spectrum} and Theorem~\ref{esv-spectrum} shows how sensitive the range of possible automorphism groups is to the choice of axioms, even among models satisfying Foundation. With Infinity present, the role of Choice in this connection was already apparent in early investigations of automorphisms of models of set theory. Cohen \cite[Theorem 1]{Cohen} constructed models of \textup{\textsf{ZF}} admitting nontrivial automorphisms of order two. In his doctoral work under Keisler, announced in 1969, five years before the publication of Cohen's paper, Morris \cite{Morris,MorrisThesis} showed that every countable group embeds into the automorphism group of a suitable model of $\ZF$. In his announcement \cite{Morris}, Morris notes that the special case of involutions answers a question posed by Cohen.\footnote{I am grateful to Asaf Karagila for bringing Morris's work on automorphisms of models of set theory to my attention. Morris's dissertation also contains a different construction, now known as the Morris model, of a model of $\ZF$ that cannot be extended to a model of $\mathsf{ZFC}$ without adding ordinals. See Karagila \cite{Karagila} for a modern treatment of this construction and an account of Morris's original approach.} An alternative, syntactic proof of the existence of models of $\ZF$ admitting involutions was given by Anapolitanos \cite{An1}.

In the presence of Choice, however, the possibilities are more restricted. As noted by Friedman (see \cite[p.~326]{Cohen}), every automorphism of a model of \textup{\textsf{ZFC}} is uniquely determined by its action on the ordinals; in other words, the action of the automorphism group on the ordinals is \emph{faithful}. Faithfulness also holds in the presence of the Leibniz--Mycielski axiom (\textup{\textsf{LM}})\footnote{The Leibniz--Mycielski axiom, introduced in
\cite{LM}, is a set-theoretic formulation of Leibniz's
principle of the identity of indiscernibles. It asserts
that any two distinct sets can be distinguished by a
parameter-free formula in some rank segment containing
both of them. By \cite[Theorems~3.1 and~3.3]{LM} (the latter due to Solovay), over $\ZF$, $\mathsf{LM}$ and $\mathsf{AC}$
are mutually independent, but both are implied by
$\mathsf{V=OD}$. See the proof of
Proposition~\ref{faithful} for an equivalent formulation
of $\mathsf{LM}$ useful in the study of automorphisms.} \cite[proof of Theorem 3.6(2)]{Paris}. Both assertions are stated and proved in Proposition~\ref{faithful}. Since the ordinals are linearly ordered, faithfulness entails right-orderability (see Theorem~\ref{orderability}). What remains to be determined is whether every right-orderable group can be realized as the full automorphism group of a model of each of $\mathsf{ZFC}$ and $\ZF+\mathsf{LM}$.

Our first result answers this question affirmatively. In fact, realization can be accomplished in models of $\ZF+\mathsf{V=OD}$, where $\mathsf{V=OD}$ asserts that every set is ordinal definable. This assertion yields a parameter-free definable well-ordering of the universe, and hence the availability of definable Skolem functions (Definition~\ref{def:skolem}). Moreover, over $\ZF$, $\mathsf{V=OD}$ implies both $\mathsf{AC}$ and $\mathsf{LM}$ (see Subsection~\ref{sec:applications}).

The spectrum question concerns the class of automorphism groups associated with a theory, without specifying the complete theory of the models witnessing realization. Can we prescribe the group while retaining the elementary theory of a given model? In view of Gaifman's theorem and Theorem~\ref{ks}, it is also natural to ask whether some form of the end-extension property can be retained. 

There is, however, a serious obstacle to obtaining a literal set-theoretic counterpart of the end-extension results for arithmetic. Whereas the MacDowell--Specker theorem \cite[Theorem 2.2.8]{KS} guarantees that every model of \textup{\textsf{PA}} has a proper elementary end extension, Kaufmann \cite[Theorem~4.2]{Kaufmann} showed that every consistent extension of \textup{\textsf{ZFC}} in the language of set theory has a model of cardinality $\aleph_1$ with no proper elementary end extension (see also \cite[Section 2]{Bogota}).\footnote{More recently, the author \cite[Theorem 5.18]{Dead} showed that every countable model of $\ZF$ has an elementary extension of cardinality $\aleph_1$ with no proper end extension to a model of $\ZF$.} In part (2) of Theorem~A below, the elementary end embedding concerns the canonically interpreted arithmetic structures and preserves the restrictions of all parameter-free definable relations of the ambient models. Via Ackermann coding, it induces an elementary end embedding between the corresponding hereditarily finite structures.

For a model $\mathcal K$ of $\ZF$, let $\mathcal A_{\mathcal K}^{+}$ denote its canonically interpreted arithmetic structure expanded by predicates for the restrictions of all parameter-free definable relations of $\mathcal K$ (Definition~\ref{def:kunen}).

\begin{theoremA}
\begin{enumerate}[label=(\arabic*)]
\item The following equalities and strict inclusion hold:
\[
\begin{aligned}
\Grpr&=\G(\PA)=\G(\ZF+\mathsf{V=OD})\\
&=\G(\mathsf{ZFC})=\G(\ZF+\mathsf{LM})\\
&\subsetneq\G(\ZF)\subseteq\Grp.
\end{aligned}
\]
\item Suppose $\mathcal M\models\ZF+\mathsf{V=OD}$ and $G\in\Grpr$. There is a model $\mathcal N\equiv\mathcal M$ with $\Aut(\mathcal N)\cong G$ and an elementary end embedding
\[
j:\mathcal A_{\mathcal M}^{+}\longrightarrow\mathcal A_{\mathcal N}^{+}
\]
which, via Ackermann coding, induces an elementary end embedding
\[
j_{\mathrm{HF}}:(V_\omega^{\mathcal M},\in^{\mathcal M})\longrightarrow(V_\omega^{\mathcal N},\in^{\mathcal N}).
\]
\end{enumerate}
\end{theoremA}

Corollaries~\ref{od-spectrum} and~\ref{choice-spectrum} establish part (1). Theorem~\ref{hfend} gives a strengthened version of part (2), including Skolem generation and the unique extension of the arithmetic embedding to the Skolem hull of the original arithmetic domain. 

The proof of Theorem A rests on two complementary
constructions. The first enriches a designated domain
of a parent model by naming the restrictions of its
parameter-free definable relations. We call the resulting
structure the \emph{induced expansion} (Definition~\ref{def:kunen}).
Our approach is motivated by Kunen's suggestion acknowledged
by Abramson and Harrington \cite[p.~589, footnote~9]{AH},
concerning the extension of their results on models without
indiscernibles from arithmetic to theories with definable
Skolem functions in which arithmetic can be interpreted
with induction for all predicates definable in the parent
model. In this arithmetic setting we refer to the induced
expansions as \emph{Kunen expansions}.
In our set-theoretic applications, the parent model
$\mathcal M$ is a model of $\ZF$, and the designated domain
is $\omega^{\mathcal M}$, the natural numbers in the sense of
$\mathcal M$. The corresponding Kunen expansion, with the
canonical arithmetic operations included, is precisely the
structure $\mathcal A_{\mathcal M}^{+}$ appearing in Theorem A(2).

The second construction concerns parent models that
are \emph{definably generated} by their designated
domains, i.e., every element of the parent model is
uniquely definable from a finite tuple of elements
of its designated domain (Definition~\ref{def:definable-generation}). For a theory with definable
Skolem functions, the construction associates with
a model $\mathcal B^+$ of the theory of induced
expansions a term model $\mathcal K(\mathcal B^+)$,
which we call the \emph{reconstructed model}
(the \emph{Kunen model} in the arithmetic setting described above). Its elements are
equivalence classes of formal Skolem terms, with equality and
relations determined by the named predicates of $\mathcal B^+$.
When $\mathcal B^+$ arises from a parent model, the reconstructed
model is canonically isomorphic to the Skolem hull of the
designated domain:
\[
\mathcal K(\mathcal P_{\mathcal M}^{+})\cong\Sk^{\mathcal M}(P^{\mathcal M}).
\]
Here $\Sk^{\mathcal M}(B)$ denotes the Skolem hull
of $B$ in $\mathcal M$.
In particular, for $\mathcal M\models\ZF+\mathsf{V=OD}$,
we have
\[
\mathcal K(\mathcal A_{\mathcal M}^{+})
\cong\Sk^{\mathcal M}(\omega^{\mathcal M}).
\]

The existence of a parent model for a prescribed
countable induced expansion follows from the omitting
types theorem (Proposition~\ref{kunen-countable-existence}),
without any assumption concerning definable Skolem
functions. The alternative proof of Theorem~4.2 of
Abramson and Harrington \cite[p.~590]{AH} contains the
arithmetic special case of our Reconstruction Theorem
(Theorem~\ref{kunen-reconstruction}): definable Skolem
functions yield definably generated parent models for
induced expansions of arbitrary cardinality.
Theorem~\ref{transfer} shows that definable generation
alone suffices for the transfer of elementary embeddings
and automorphisms between induced expansions and their
parent models. Whereas Abramson and Harrington transfer
results concerning the absence of indiscernibles from
arithmetic to more general theories, we transfer results
concerning full automorphism groups.

The link between reconstruction and arithmetic
realization is provided by the notion of accommodation
introduced in \cite[Definition~61]{EL}. Under the hypotheses of our
general realization theorems, accommodation ensures
that Kunen expansions satisfy induction in their
expanded language, thereby making the Kossak--Schmerl
realization theorem (Theorem~\ref{ks}) applicable.
Reconstruction and transfer then enable us to lift
the arithmetic realizations to models of the parent
theory.
In this way we obtain realization theorems for
theories considerably more general than the set
theories considered in Theorem A.

The proof of Theorem B is based on a different approach suggested by the proof of Theorem~\ref{esv-spectrum}: given a group $G$, we first choose a graph $\mathcal G$ whose automorphism group is isomorphic to $G$; this is possible by classical results (see Theorem~\ref{graph-representation}). We then encode $\mathcal G$ in an extensional digraph and construct a model of set theory in which a copy of $\mathcal G$ is parameter-free definable and every automorphism of that copy extends uniquely to an automorphism of the model. Under the assumption that there is a strongly inaccessible
cardinal exceeding $|G|$, the proof of
Theorem~\ref{esv-spectrum} can be adapted to realize $G$
as the full automorphism group of a model of $\mathsf{ZFC}$
without Foundation.\footnote{When Schmerl, Visser and I wrote \cite{ESV}, we were unaware that a transfinite version of the completion construction had already appeared in the work of Tzouvaras \cite[Section~4]{Tz}. In recognition of his work, we use the term \emph{Tzouvaras completion}.} For $T=\ZF$ or $T=\mathsf{ZFC}$, we write $T^{-\mathsf{Found}}$ for the theory obtained by replacing Foundation with its negation.

\begin{theoremB}
If $\kappa$ is strongly inaccessible and $G\in\Grp$ with $|G|<\kappa$, there is a model $\mathcal N$ of $\mathsf{ZFC}^{-\mathsf{Found}}$ whose externally well-founded part is isomorphic to $(V_\kappa,\in)$ and such that
\[
|N|=\kappa\quad\text{and}\quad\Aut(\mathcal N)\cong G.
\]
\end{theoremB}

The proof of Theorem B (Theorem~\ref{inacc}) proceeds by successively adjoining
new vertices to an extensional digraph to represent missing
sets, until the resulting structure is a model of
$\mathsf{ZFC}^{-\mathsf{Found}}$. Every automorphism of the
original digraph extends uniquely to an automorphism of
its completion. The crucial point is to arrange that the
original digraph is parameter-free definable in the
resulting model, so that restriction and extension induce
mutually inverse isomorphisms between their automorphism
groups. We use the coding digraph of Enayat, Schmerl and Visser \cite[Theorem 4.2]{ESV}, described in Subsection~\ref{subsec:recovery}, and establish its definability in the completed structure in Lemma~\ref{seed}. 

In Appendix~\ref{sec:class-completions} we explain
how a model of $\ZF$ and a prescribed parametrically
definable extensional digraph whose fibres are sets
in that model can both be end embedded into a common
model of set theory without Foundation. We use this result to
revisit a question concerning external descending
membership chains in models with well-founded ordinals,
previously addressed by the author in \cite{MO}.

The paper is organized as follows.
\begin{itemize}
\item Section~\ref{sec:prelim} establishes notation and conventions and presents the representation theorems for ordered structures and graphs.
\item Section~\ref{sec:arithmetic-transfer} develops
the theory of induced expansions and reconstructed models,
establishing existence, reconstruction and transfer
results, and examining their implications for
bi-interpretability.
\item Section~\ref{sec:realization} combines these results with arithmetic realization to establish Theorem A.
\item Section~\ref{sec:nofoundation} establishes Theorem B by Tzouvaras completion.
\item Section~\ref{sec:open-question} poses open questions concerning models of $\ZF$ and its extensions.
\item Appendix~\ref{sec:class-completions} examines proper-class completions and descending membership chains.
\end{itemize}

\section{Preliminaries}
\label{sec:prelim}
We establish notation for the structures interpreted in models of set theory and present the group representation results used in Sections~\ref{sec:realization} and~\ref{sec:nofoundation}.

\subsection{Notation and conventions}\label{subsec:notation}
We work in $\mathsf{ZFC}$ in the metatheory unless otherwise indicated. Unless specified otherwise, theories under consideration are consistent, and finiteness and cardinality are understood externally.

We use italic letters for theory variables and upright sans serif for names of specific theories. For a structure $\mathcal M$, we write $\operatorname{Th}(\mathcal M)$ for its complete theory.

We write $\Lset=\{\in\}$ for the usual language of set theory (with equality).

We use calligraphic letters for structures and the corresponding Roman letters for their universes. For an $\Lset$-structure $\mathcal M=(M,E)$ and $a\in M$, the \emph{fibre} of $a$ is the set of its members in the sense of $\mathcal M$:
\[
 a_E=\{x\in M:xEa\}.
\]
For $\mathcal M\models\ZF\setminus\{\mathsf{Found}\}$, we write $\operatorname{Ord}^{\mathcal M}$, $\omega^{\mathcal M}$ and $V_\omega^{\mathcal M}$ for the ordinals, natural numbers and hereditarily finite sets, respectively, in the sense of $\mathcal M$. We denote by $\mathcal A_{\mathcal M}$ the arithmetic structure with universe $\omega^{\mathcal M}$ and the canonically interpreted arithmetic operations, constants and order. The corresponding hereditarily finite structure is $(V_\omega^{\mathcal M},\in^{\mathcal M})$, with $\in^{\mathcal M}$ restricted to the indicated universe. The expansion $\mathcal A_{\mathcal M}^{+}$ used in Theorem A(2) is a special case of Definition~\ref{def:kunen}.

We write $\PA^{*}$ for Peano arithmetic in a countable expansion of the language of arithmetic, with induction for all formulae of the expanded language. We typically denote models of $\PA^{*}$ by $\mathcal A^{+}$, with $\mathcal A$ denoting their arithmetic reducts and $A$ their common universe.

We say \emph{$\mathcal M$-definable} for definability without parameters, and \emph{parametrically $\mathcal M$-definable} when parameters are allowed.

\begin{definition}\label{def:standardness}
\begin{enumerate}[label=(\alph*)]
\item An $\Lset$-structure $\mathcal M$ is called \emph{standard} if $\in^{\mathcal M}$ is externally well-founded.
\item A model $\mathcal M$ of $\ZF\setminus\{\mathsf{Found}\}$ is called \emph{ordinal-standard} if $(\operatorname{Ord}^{\mathcal M},\in^{\mathcal M})$ is externally well-founded, in which case its order type is called the \emph{ordinal height} of $\mathcal M$.
\end{enumerate}
\end{definition}

For models of $\ZF$, standardness and ordinal-standardness coincide, by the availability of the rank function.

\Needspace{6\baselineskip}
\begin{definition}\label{def:end-extension}
Let $\mathcal M\subseteq\mathcal N$.
\begin{enumerate}[label=(\alph*)]
\item For $\Lset$-structures $\mathcal M=(M,E)$ and $\mathcal N=(N,F)$, $\mathcal N$ is an \emph{end extension} of $\mathcal M$ if
\[
\forall a\in M\quad a_E=a_F.
\]
\item For arithmetic structures, $\mathcal N$ is an \emph{end extension} of $\mathcal M$ if
\[
\forall a\in M\ \forall b\in N\ (b<^{\mathcal N}a\Rightarrow b\in M).
\]
\end{enumerate}
In either setting, an \emph{end embedding} is an embedding whose image is an end submodel. We write $\prec_{\mathrm{end}}$ for elementary end extension.
\end{definition}

\subsection{Automorphism groups: topology and representation theorems}\label{subsec:topology}\label{subsec:auttop}
The transfer and realization theorems of Sections~\ref{sec:arithmetic-transfer} and~\ref{sec:realization} have topological as well as abstract group-theoretic conclusions. We endow $\Aut(\mathcal M)$ with the topology of pointwise convergence on the discrete universe $M$. For $B\subseteq M$, the \emph{pointwise stabilizer} of $B$ is
\[
\Aut(\mathcal M/B)=\{g\in\Aut(\mathcal M):g\restriction B=\operatorname{id}_B\}.
\]
The subgroups $\Aut(\mathcal M/F)$, for finite $F\subseteq M$, form a neighborhood basis at the identity. Equivalently, this is the topology inherited from the product space $M^M$; the automorphism group is closed in the symmetric group on $M$. We write $\cong_{\mathrm{top}}$ for topological group isomorphism and $\cong$ for abstract group isomorphism when comparing groups.

The connection between right-orderability and automorphisms of linearly ordered structures rests on the following classical characterization.
\begin{theorem}[{Cohn--Conrad \cite{Cohn,Conrad};
see also \cite[p.~1249]{Schmerl}}]
\label{orderability}
For a group $G$, the following conditions are equivalent:
\begin{enumerate}[label=(\arabic*)]
\item $G$ is right-orderable.
\item $G$ is left-orderable.
\item $G$ is isomorphic to a subgroup of the
automorphism group of a linear order.
\item $G$ is isomorphic to a closed subgroup of the
automorphism group of a linear order.
\item $G$ is isomorphic to the automorphism group of
a linearly ordered structure $(K,<,R)$ with one
additional binary relation $R$ and $|K|=|G|$.
\end{enumerate}
\end{theorem}
In particular, right-orderable groups are torsion-free.

In contrast to linear orders, graphs impose no restriction on their full automorphism groups. The following theorem provides the starting point for the construction in Section~\ref{sec:nofoundation}.
\begin{theorem}[Frucht, de Groot, Sabidussi, Babai]\label{graph-representation}
Every group $G$ is isomorphic to the automorphism group of a simple undirected graph $\mathcal A=(A,F)$. Moreover, $\mathcal A$ can be chosen finite if $G$ is finite, and with $|A|=|G|$ if $G$ is infinite.\footnote{Frucht \cite{Frucht} established the finite case, answering a question of K\"onig. De Groot \cite{deGroot} and Sabidussi \cite{Sabidussi} independently extended the result to arbitrary groups. The cardinality assertion for infinite groups follows from Babai's stronger result \cite[p.~26, Corollary]{Babai}, which gives a simple graph with three vertex orbits.}
\end{theorem}

\section{Induced expansions: existence, reconstruction and transfer}\label{sec:arithmetic-transfer}

The induced expansion construction relates a definable domain, enriched by
the relations it inherits from its parent model, to the
structure of the parent model itself. We develop the
structural properties of this construction independently of
the arithmetic realization theorems, before combining the two
in Section~\ref{sec:realization}. Proposition~\ref{kunen-countable-existence}
separates countable existence from the stronger reconstruction
results that use definable Skolem functions.
The applications in Section~\ref{sec:realization}
rely on Propositions~\ref{accommodation-induction} and~\ref{kunen-internal-reconstruction}
and Theorems~\ref{kunen-reconstruction} and~\ref{transfer},
whereas the uniqueness, categorical and infinitary bi-interpretability
results (Corollary~\ref{kunen-categorical} and
Remarks~\ref{kunen-functor} and~\ref{kunen-biinterpretability})
are of independent interest and are not invoked
in the subsequent sections.

\subsection{Accommodation and induced expansions}
\label{subsec:accommodation}\label{sec:transfersec}

We first consider expansions induced on arbitrary definable
domains. When the domain carries an interpretation of
arithmetic, accommodation ensures that the resulting
expansion satisfies full induction.

\begin{definition}\label{def:kunen}
Let $T$ be a theory in a countable language $\mathcal L$,
and let $P(x)$ be a parameter-free formula such that
\[
T\vdash\exists x\,P(x).
\]
For each $\mathcal M\models T$, we call $\mathcal M$ the
\emph{parent model} and refer to $P^{\mathcal M}$ as the
\emph{designated domain}.

The \emph{induced expansion} $\mathcal P_{\mathcal M}^{+}$
has universe $P^{\mathcal M}$ and, for each
$\mathcal L$-formula $\varphi(\bar x)$, a predicate
$R_\varphi$ interpreted by
\[
R_\varphi^{\mathcal P_{\mathcal M}^{+}}
=
\{\bar a\in(P^{\mathcal M})^{|\bar x|}:
\mathcal M\models\varphi(\bar a)\}.
\]
We allow nullary predicates corresponding to sentences
of $\mathcal L$.

If $P$ is the domain of a parameter-free, one-dimensional,
equality-preserving interpretation $\sigma$ of arithmetic,
we call $\mathcal A=\sigma(\mathcal M)$ the
\emph{child model} of $\mathcal M$. We denote its induced
expansion, including the interpreted arithmetic operations,
by $\mathcal A_{\mathcal M}^{+}$. When $\sigma$ witnesses
accommodation of $\PA$, we also call this expansion the
\emph{Kunen expansion}.
\end{definition}

The arithmetic terminology recognizes Kunen's suggestion
acknowledged in \cite[p.~589, footnote~9]{AH}. The general
construction belongs to the framework of model theory over a
predicate and can be viewed
as the restriction of the Morleyization of the parent model
to its designated domain; see \cite[pp.~62--64]{Hodges}.

Let $\mathcal L_P^+$ denote the language of the induced expansions, including the interpreted arithmetic symbols in the arithmetic case. In statements applying to both cases, $\mathcal P_{\mathcal M}^{+}$ denotes the corresponding expansion. For each $\mathcal L_P^+$-formula $\theta$,
we define its translation $\theta^\dagger$ into
$\mathcal L$ recursively. For atomic formulas with variables as arguments in the general relational case, and for the logical connectives and quantifier, the clauses are as follows:
\begin{align*}
(x=y)^\dagger&=(x=y),\\
(R_\varphi(\bar x))^\dagger&=\varphi(\bar x),\\
(\neg\theta)^\dagger&=\neg\theta^\dagger,\\
(\theta\land\psi)^\dagger
&=\theta^\dagger\land\psi^\dagger,\\
(\exists x\,\theta)^\dagger
&=\exists x(P(x)\land\theta^\dagger).
\end{align*}
Bound variables are renamed when necessary to avoid unintended clashes.

In the arithmetic case, for each arithmetic term
$t(\bar x)$, let $\delta_t(z,\bar x)$ define its
interpreted value in the parent model. Thus, for example,
\begin{align*}
\delta_{x_i}(z,\bar x)&\equiv z=x_i,\\
\delta_{s+t}(z,\bar x)&\equiv
\exists u\exists v\bigl(
\delta_s(u,\bar x)\land
\delta_t(v,\bar x)\land
\delta_+(z,u,v)\bigr),
\end{align*}
where $\delta_+$ defines the graph of the interpreted
addition operation. The remaining arithmetic operations
and constants are treated analogously.

For an atomic formula $Q(t_1,\ldots,t_k)$, let $q(\bar z)$ be $\varphi(\bar z)$ if $Q=R_\varphi$, the formula $z_1=z_2$ if $Q$ is equality, and the ambient defining formula of $Q$ if $Q$ is an arithmetic relation symbol. Put
\[
\begin{aligned}
(Q(t_1,\ldots,t_k))^\dagger
\equiv\exists z_1\cdots\exists z_k
\Bigl(&\bigwedge_{i=1}^{k}\delta_{t_i}(z_i,\bar x)\\
&\land q(z_1,\ldots,z_k)\Bigr).
\end{aligned}
\]

\begin{proposition}\label{kunen-translation}
For every parent model $\mathcal M\models T$, every
formula $\theta(\bar x)$ of the language of its induced expansion, and every tuple $\bar a$ from its designated
domain,
\begin{equation}\label{eq:kunen-translation}
\mathcal P_{\mathcal M}^{+}\models\theta(\bar a)
\quad\Longleftrightarrow\quad
\mathcal M\models\theta^\dagger(\bar a).
\end{equation}
\end{proposition}

\begin{proof}
We use induction on the complexity of $\theta$. The atomic cases
follow from the definitions, and the Boolean cases are
immediate. For the existential quantifier, writing
$B=P^{\mathcal M}$, we have
\begin{align*}
\mathcal P_{\mathcal M}^{+}
\models\exists x\,\theta(x,\bar a)
&\iff
\exists b\in B\
\mathcal P_{\mathcal M}^{+}\models\theta(b,\bar a)\\
&\iff
\exists b\in B\
\mathcal M\models\theta^\dagger(b,\bar a)\\
&\iff
\mathcal M\models
\exists x(P(x)\land\theta^\dagger(x,\bar a)).
\end{align*}

\end{proof}

Let $S_T$ be the common theory of the induced expansions
of models of $T$, with the designated domain fixed.
By \eqref{eq:kunen-translation} and the completeness
theorem,
\begin{equation}\label{eq:kunen-theory}
S_T=
\{\theta\in\operatorname{Sent}(\mathcal L_P^+):
T\vdash\theta^\dagger\}.
\end{equation}
In particular, if $T$ is complete, then $S_T$ is complete,
and for any $\mathcal M\models T$,
\[
S_T=\operatorname{Th}(\mathcal P_{\mathcal M}^{+}).
\]

\begin{corollary}\label{kunen-qe}
The theory $S_T$ admits elimination of quantifiers.
\end{corollary}

\begin{proof}
For every $\mathcal M\models T$ and every tuple $\bar a$
from its designated domain,
\begin{align*}
\mathcal P_{\mathcal M}^{+}\models\theta(\bar a)
&\iff\mathcal M\models\theta^\dagger(\bar a)\\
&\iff\mathcal P_{\mathcal M}^{+}
\models R_{\theta^\dagger}(\bar a).
\end{align*}
Consequently,
\[
S_T\vdash\forall\bar x
\bigl(\theta(\bar x)\leftrightarrow
R_{\theta^\dagger}(\bar x)\bigr).
\]
Thus every formula is equivalent modulo $S_T$ to an
atomic formula.
\end{proof}

We now specialize to interpretations of arithmetic.
Accommodation strengthens interpretability by requiring
induction not only for arithmetically definable predicates
but also for predicates inherited from the parent model.

\begin{definition}\label{def:accommodation}
Suppose $\mathcal M$ interprets a model $\mathcal A$ of $\PA$
via an interpretation $\sigma$, possibly with parameters.
If $P\subseteq M^k$ is the interpretation domain and $E$
is the equivalence relation interpreting equality, then $A=P/E$.
We say that $\mathcal M$ \emph{accommodates} $\PA$ via $\sigma$
if the following condition holds.

For every ambient formula $\varphi(\bar x,\bar y)$
and every tuple $\bar b$ from $M$ for which
$\varphi(\bar x,\bar b)$ is invariant under $E$ on $P$,
put
\[
X_{\varphi,\bar b}
=
\{[\bar a]_E:\bar a\in P,\
\mathcal M\models\varphi(\bar a,\bar b)\}.
\]
Then
\[
\begin{aligned}
0^{\mathcal A}\in X_{\varphi,\bar b}
\ \land\
S^{\mathcal A}[X_{\varphi,\bar b}]
\subseteq X_{\varphi,\bar b}
\quad\Longrightarrow\quad
X_{\varphi,\bar b}=A.
\end{aligned}
\]
Equivalently, the child model has no proper cut induced
by a parametrically definable predicate of the parent model.

A theory $T$ \emph{accommodates} $\PA$ via an interpretation
$\sigma$ of $\PA$ in $T$ if every model of $T$ accommodates
$\PA$ via $\sigma$.\footnote{Accommodation is defined for arbitrary
scheme templates in \cite[Definition 61]{EL}. The
embedding of $\PA$ in the sense of
\cite[pp.~589--590]{AH} is the special case in which
the interpretation is parameter-free, one-dimensional
and equality-preserving.}
Here a theory interpretation is uniform and parameter-free;
a uniform interpretation with parameters can also be used
by specifying a parameter condition provably realized in $T$
and requiring accommodation for every tuple satisfying that condition.
\end{definition}

A parametrically definable relation on $P^n$ induces a
relation on $A^n$ if it is invariant under $E$ in each
argument. We have the following characterization of
accommodation for arbitrary interpretations, whose proof is straightforward. 

\begin{proposition}\label{accommodation-induction}
Suppose $\mathcal M$ parametrically interprets a model
$\mathcal A$ of $\PA$. The following conditions are equivalent:
\begin{enumerate}[label=(\alph*)]
\item $\mathcal M$ accommodates $\PA$ via the given interpretation.
\item $\mathcal A$ satisfies induction in every expansion
by finitely many relations induced by parametrically
definable relations of $\mathcal M$.
\item $\mathcal A$ satisfies full induction in its expansion
by all relations induced by parametrically definable
relations of $\mathcal M$.
\end{enumerate}
\end{proposition}

In particular, if $T$ accommodates $\PA$ via a parameter-free,
one-dimensional, equality-preserving interpretation $\sigma$,
then
\[
\mathcal A_{\mathcal M}^{+}\models\PA^*
\qquad\text{for every }\mathcal M\models T.
\]
Consequently, by the completeness theorem,
\[
S_T\supseteq\PA^*.
\]
This is the property of accommodation that enables us
to apply Schmerl's realization theorems to Kunen expansions in Section~\ref{sec:realization}.

Although accommodation does not require interpretations
to be one-dimensional or equality-preserving, the reconstruction
in the next subsection uses designated domains consisting of
individual elements of the parent model, which serve as arguments
of its definable Skolem functions. The following proposition
identifies circumstances under which interpretations can be
put in this form.

Here \emph{elimination of imaginaries} means that every parameter-free
definable equivalence relation on a parameter-free definable
subset of a finite Cartesian power of the universe admits a
parameter-free definable coding function into the universe,
with equal values precisely on equivalent tuples. The defining
formula is required to work uniformly in all models of $T$.

\begin{proposition}\label{interpretation-reductions}
Suppose $T$ accommodates $\PA$ via a parameter-free interpretation
$\sigma$.
\begin{enumerate}[label=(\alph*)]
\item If $T$ has a parameter-free definable pairing
function, then $\sigma$ can be replaced by a
one-dimensional interpretation witnessing accommodation.
\item If $T$ eliminates imaginaries, then $\sigma$
can be replaced by a one-dimensional, equality-preserving
interpretation witnessing accommodation.
\end{enumerate}
\end{proposition}

\begin{proof}
Work in an arbitrary parent model $\mathcal M\models T$.
Let $P\subseteq M^k$ be the interpretation domain and
$E$ the equivalence relation interpreting equality.

(a) Iteration of the pairing function gives a definable
injection
\[
p_k:M^k\longrightarrow M.
\]
Put $P'=p_k[P]$ and define $E'$ on $P'$ by
\[
p_k(\bar a)\ E'\ p_k(\bar b)
\quad\Longleftrightarrow\quad
\bar a\ E\ \bar b.
\]
Transporting the interpreted operations and relations
along $p_k$ gives an interpretation whose child model
is canonically isomorphic to the original child model.

(b) By elimination of imaginaries, there is a
parameter-free definable function $q$ on $P$ such that
\[
q(\bar a)=q(\bar b)
\quad\Longleftrightarrow\quad
\bar a E\bar b.
\]
Put $P_0=q[P]$. The map
\[
P/E\longrightarrow P_0,\qquad
[\bar a]_E\longmapsto q(\bar a)
\]
is a bijection. Transporting the interpreted relations
and functions along this bijection yields a one-dimensional,
equality-preserving interpretation.

In both cases the resulting child model is isomorphic
to the original one via an ambient-definable
correspondence. Parametrically definable predicates
pull back along this correspondence, so accommodation
is preserved.
\end{proof}

\begin{remark}
$\ZF$ accommodates $\PA$ via the canonical
parameter-free interpretation of arithmetic on $\omega$,
by Separation. Moreover, $\ZF$ satisfies both hypotheses
of Proposition~\ref{interpretation-reductions}, since
it has a parameter-free definable pairing function and
eliminates imaginaries by Scott's trick (see
\cite[Remark~49]{EL} and
\cite[Section~4.4, Exercise~2]{Hodges}).
Every subsystem of $\mathsf Z_2$ containing the arithmetic
axioms and full induction also accommodates $\PA$ via its
number sort.\footnote{More generally, all inductive
sequential theories accommodate $\PA$. Recall that a
sequential theory is called \emph{inductive} if it admits
a parameter-free interpretation of arithmetic satisfying
induction for all predicates induced by parametrically
definable relations of the parent model. By Montague's
theorem \cite{Montague}, no consistent inductive sequential
theory is finitely axiomatizable. See also
\cite[Section 3.2]{EL} and
\cite[Definition 2.1(g)]{ESV}.}
\end{remark}
\subsection{Reconstruction and transfer}
\label{sec:generalsec}\label{subsec:reconstruction}

The preceding subsection establishes the relationship
between definable relations in a parent model and those
in its induced expansion. We first address the existence of
parent models for prescribed induced expansions, and then
examine the extent to which such expansions determine their
parent models. Accommodation plays no role in these results.

\begin{definition}\label{def:definable-generation}
Let $\mathcal M$ be a structure and $B\subseteq M$.
We say that $\mathcal M$ is \emph{definably generated}
by $B$ if every element of $M$ is uniquely definable
in $\mathcal M$ from a finite tuple of elements of $B$.
\end{definition}

The following proposition is a standard application
of the omitting types theorem in the context of
classification over a predicate. In the terminology
of this literature, it establishes the Gaifman
property in countable cardinalities
(see \cite[Lemmas~1.6--1.7]{Pillay2026}).

\begin{proposition}\label{kunen-countable-existence}
Let $T$ be a theory in a countable language and $P(x)$
a parameter-free formula such that $T\vdash\exists x\,P(x)$.
Every countable model $\mathcal B^+\models S_T$ is the
induced expansion of the designated domain of some countable
model $\mathcal N\models T$, up to isomorphism.
\end{proposition}

\begin{proof}
Introduce constants $c_b$ for $b\in B$, and let $\Gamma$
consist of $T$, the sentences $P(c_b)$, and the translations
$\theta^\dagger(\bar c_{\bar b})$ of the elementary diagram
of $\mathcal B^+$. Every finite conjunction of the latter
sentences, with the constants replaced by variables and
existentially quantified over $P$, is consistent with $T$
by \eqref{eq:kunen-theory}. Thus $\Gamma$ is consistent
by compactness. Consider the partial type
\[
\Sigma(x)=\{P(x)\}\cup\{x\ne c_b:b\in B\}.
\]
If a formula $\varphi(x,\bar c_{\bar b})$ consistent with
$\Gamma$ implied every member of $\Sigma(x)$ modulo $\Gamma$,
then $\mathcal B^+\models\exists x\,R_\varphi(x,\bar b)$:
otherwise the translated elementary diagram would contradict
its consistency together with $P(x)$. Hence for some $d\in B$,
$\mathcal B^+\models R_\varphi(d,\bar b)$, so that
$\Gamma\vdash\varphi(c_d,\bar c_{\bar b})$, a contradiction.
By the omitting types theorem, $\Gamma$ has a countable model
omitting $\Sigma$. Its reduct to the language of $T$ is the
desired model $\mathcal N$, since its designated domain consists
precisely of the interpretations of the constants $c_b$.
\end{proof}

\begin{remark}
The existence problem for induced expansions belongs to the
framework of classification over a predicate, originating
in Gaifman's work \cite{Gaifman1974}; see also
\cite{PillayShelah}, \cite[Section~12.5]{Hodges}, and
\cite{Pillay2026,Usvyatsov}. Proposition~\ref{kunen-countable-existence}
establishes countable existence without definable Skolem
functions. Countable existence does not, by itself, ensure that the
parent model is definably generated by its designated domain.
The availability of definable Skolem functions enables us to
obtain such parent models in arbitrary cardinalities, as shown
in Theorem~\ref{kunen-reconstruction}. Definable generation alone suffices
for the transfer of elementary embeddings and automorphisms.
This does not assert relative categoricity for all models of $T$.
\end{remark}

\begin{definition}\label{def:skolem}
A theory $T$ in a language $\mathcal L$ has \emph{definable Skolem functions} if,
for every $\mathcal L$-formula $\varphi(x,\bar y)$, with distinguished witness variable $x$ and a specified tuple $\bar y$ of distinct parameter variables (not including $x$), containing all other free variables, there is an $\mathcal L$-formula
$\delta_\varphi(x,\bar y)$ such that
\[
\begin{aligned}
T&\vdash\forall\bar y\exists!x\,
\delta_\varphi(x,\bar y),\\
T&\vdash\forall x\forall\bar y
\bigl[
((\exists z\,\varphi(z,\bar y))\land
\delta_\varphi(x,\bar y))
\longrightarrow\varphi(x,\bar y)
\bigr].
\end{aligned}
\]
For each such indexed formula $\varphi$, fix once and for all a formula $\delta_\varphi$
as above and introduce a fresh function symbol
$f_\varphi$ of arity $|\bar y|$.
Every model $\mathcal M$ of $T$ has a unique
expansion interpreting $f_\varphi$ as the total
function $f_\varphi^{\mathcal M}$ given by
\[
f_\varphi^{\mathcal M}(\bar a)=b
\quad\Longleftrightarrow\quad
\mathcal M\models\delta_\varphi(b,\bar a).
\]
\end{definition}

Suppose $T$ has definable Skolem functions,
$\mathcal M\models T$, and $B\subseteq M$.
The \emph{Skolem hull} of $B$ in $\mathcal M$,
denoted by $\Sk^{\mathcal M}(B)$, is the least
elementary submodel of $\mathcal M$ containing $B$.
Its universe consists precisely of the elements
uniquely definable in $\mathcal M$ from finite tuples
of elements of $B$. Thus
$\mathcal M=\Sk^{\mathcal M}(B)$ if and only if
$\mathcal M$ is definably generated by $B$
(Definition~\ref{def:definable-generation}).

In particular, if $T$ is complete, the Skolem hull
of the empty set in any model of $T$ is a prime
model of $T$.

The underlying idea of reconstruction is to regard the Skolem hull as a term model. Let $T$ be a theory in a countable language $\mathcal L$ with definable Skolem functions. Fix an enumeration $(\varphi_i(x,\bar y_i):i<\omega)$ of the formulas indexed as in Definition~\ref{def:skolem}, with their distinguished witness variable $x$ and specified tuples of parameter variables $\bar y_i$ (possibly empty). Put
\[
f_i:=f_{\varphi_i},\qquad
\delta_i:=\delta_{\varphi_i},\qquad
n_i:=|\bar y_i|.
\]
Thus $f_i$ is a formal function symbol of arity $n_i$, and $\delta_i$ is its defining $\mathcal L$-formula. We rename the free variables of $\delta_i$ as needed. Let $i_0$ be the index of the formula $x=y$ with parameter variable $y$. By the Skolem witness condition, $f_{i_0}^{\mathcal M}(a)=a$ for every $\mathcal M\models T$ and every $a\in M$.

Fix a parent model $\mathcal M\models T$ and let $B=P^{\mathcal M}$ be its designated domain. Put
\[
\operatorname{Term}(B):=\{\langle i,\bar a\rangle:i<\omega,\ \bar a\in B^{n_i}\}.
\]
We write $f_i(\bar a)$ for the formal term $\langle i,\bar a\rangle$, distinguishing it from its value $f_i^{\mathcal M}(\bar a)$ in the parent model. For each $i,j<\omega$, put
\[
\eta_{ij}(\bar u,\bar v):=\exists z\,(\delta_i(z,\bar u)\land\delta_j(z,\bar v)).
\]
Define an equivalence relation on terms by
\[
f_i(\bar a)\sim f_j(\bar b)\iff
\mathcal P_{\mathcal M}^{+}\models R_{\eta_{ij}}(\bar a,\bar b).
\]
For each ambient $k$-ary relation symbol $R$, put
\[
\chi^R_{\bar i}(\bar u_1,\ldots,\bar u_k):=
\exists y_1\cdots\exists y_k\left(
\bigwedge_{\ell=1}^k\delta_{i_\ell}(y_\ell,\bar u_\ell)
\land R(y_1,\ldots,y_k)\right).
\]
The relation on equivalence classes is given by
\[
R^*([f_{i_1}(\bar a_1)]_\sim,\ldots,[f_{i_k}(\bar a_k)]_\sim)
\iff \mathcal P_{\mathcal M}^{+}\models R_{\chi^R_{\bar i}}(\bar a_1,\ldots,\bar a_k).
\]
Function symbols (including constants) are handled by their graphs. Let $\mathcal K(\mathcal P_{\mathcal M}^{+})$ be the $\mathcal L$-structure with universe
$\operatorname{Term}(B)/{\sim}$ and interpretations of the symbols of
$\mathcal L$ as specified above. The evaluation map
\[
\varepsilon:\mathcal K(\mathcal P_{\mathcal M}^{+})\longrightarrow\Sk^{\mathcal M}(B),
\]
given by $[f_i(\bar a)]_\sim\longmapsto f_i^{\mathcal M}(\bar a)$, is an isomorphism of $\mathcal L$-structures by the definitions of the named predicates and the fact that every element of the Skolem hull is represented by a formal Skolem term (Definition~\ref{def:skolem}).

The construction of $\mathcal K(\mathcal P_{\mathcal M}^{+})$ depends only on the induced expansion, with $T$, $P$ and the enumeration fixed. The same formulas define $\mathcal K(\mathcal B^+)$ for any $\mathcal B^+\models S_T$: the reconstruction theorem below ensures that these definitions yield a structure. 

\begin{definition}\label{def:reconstructed-model}
We call $\mathcal K(\mathcal B^+)$ the \emph{reconstructed model} associated with $\mathcal B^+$. In the setting of accommodated arithmetic we also call it the \emph{Kunen model}.
\end{definition}

We write $\iota_B(b)=[f_{i_0}(b)]_\sim$ for the canonical identification of $B$ with its image in the reconstructed model.

This external construction requires neither arithmetic on $B$ nor passage to $S_T$. Remark~\ref{kunen-biinterpretability} explains how arithmetic coding turns it into an infinitary interpretation. The following proposition states the key properties of the hull when the parent model is already given.

\begin{proposition}\label{kunen-internal-reconstruction}
Suppose $T$ has definable Skolem functions,
$\mathcal M\models T$, and $B=P^{\mathcal M}$ is
its designated domain. Put
\[
\mathcal M^*=\Sk^{\mathcal M}(B).
\]
Then:
\begin{enumerate}[label=(\alph*)]
\item $\mathcal M^*$ is the least elementary submodel
of $\mathcal M$ containing $B$.
\item $P^{\mathcal M^*}=B$.
\item
$\mathcal P_{\mathcal M^*}^{+}
=\mathcal P_{\mathcal M}^{+}$.
\item Evaluation gives a canonical isomorphism
\[
\mathcal K(\mathcal P_{\mathcal M}^{+})\cong\mathcal M^*.
\]
\end{enumerate}
\end{proposition}

\begin{proof}
By the definition of the Skolem hull, $\mathcal M^*\prec\mathcal M$.
If $\mathcal H\prec\mathcal M$ and $B\subseteq H$, then every
 element uniquely definable in $\mathcal M$ from a finite tuple
from $B$ belongs to $H$. Hence
\[
\Sk^{\mathcal M}(B)\prec\mathcal H,
\]
establishing (a).

By elementarity,
\[
P^{\mathcal M^*}
=M^*\cap P^{\mathcal M}
=M^*\cap B=B,
\]
which establishes (b).

For each ambient formula $\varphi(\bar x)$ and tuple
$\bar b$ from $B$,
\begin{align*}
\mathcal P_{\mathcal M^*}^{+}
\models R_\varphi(\bar b)
&\iff\mathcal M^*\models\varphi(\bar b)\\
&\iff\mathcal M\models\varphi(\bar b)\\
&\iff\mathcal P_{\mathcal M}^{+}
\models R_\varphi(\bar b).
\end{align*}
In the arithmetic case, the interpreted operations
and constants agree by elementarity applied to their
defining graphs. This establishes (c). Part (d) follows from the preceding term-model construction.
\end{proof}

In particular, when the designated domain is that of
an interpreted arithmetic structure, the passage from
$\mathcal M$ to $\mathcal M^*$ leaves both the child
model and its induced expansion unchanged.

In the set-theoretic setting, the canonical Ackermann bijection is parameter-free definable in the parent model $\mathcal M$. Consequently, if $\mathcal M\models\ZF+\mathsf{V=OD}$, then
\[
\Sk^{\mathcal M}(\omega^{\mathcal M})=\Sk^{\mathcal M}(V_\omega^{\mathcal M}).
\]
Thus either designated domain can be used to reconstruct the same elementary submodel of $\mathcal M$. We work with the arithmetic domain in Section~\ref{sec:realization} in order to apply the arithmetic realization theorems directly.

Starting with an arbitrary model of $S_T$, we now obtain a
suitable parent model by compactness and pass to its Skolem hull.

\begin{theorem}[Reconstruction]
\label{kunen-reconstruction}
Suppose $T$ is in a countable language, has definable
Skolem functions, and $P(x)$ is a parameter-free formula
such that
\[
T\vdash\exists x\,P(x).
\]
For every $\mathcal B^+\models S_T$, there exists
$\mathcal N\models T$, canonically isomorphic to the reconstructed model $\mathcal K(\mathcal B^+)$, such that:
\begin{enumerate}[label=(\alph*)]
\item $P^{\mathcal N}=B$ and
$\mathcal P_{\mathcal N}^{+}=\mathcal B^+$.
\medskip
\item $\mathcal N=\Sk^{\mathcal N}(B)$.
\medskip
\item $|N|\leq\max(\aleph_0,|B|)$, with equality
$|N|=|B|$ whenever $B$ is infinite.
\end{enumerate}
\end{theorem}

\begin{proof}
Expand the language of $T$ by constants $c_b$ for
$b\in B$. Using the translation established in
Proposition~\ref{kunen-translation}, put
\[
\begin{aligned}
\Gamma=T
&\cup\{P(c_b):b\in B\}\\
&\cup\{\theta^\dagger(\bar c_{\bar b}):
\theta(\bar c_{\bar b})
\in\operatorname{ElDiag}(\mathcal B^+)\},
\end{aligned}
\]
where $\operatorname{ElDiag}(\mathcal B^+)$ denotes
the elementary diagram of $\mathcal B^+$.

Suppose $\Gamma$ is inconsistent. By compactness,
there is a formula $\theta(\bar x)$ and a tuple
$\bar b$ from $B$ such that
\[
\mathcal B^+\models\theta(\bar b)
\]
but
\[
T\vdash
\neg\exists\bar x
\left(
\bigwedge_i P(x_i)
\land\theta^\dagger(\bar x)
\right).
\]
By the definition of $S_T$ and
\eqref{eq:kunen-translation},
\[
S_T\vdash\neg\exists\bar x\,\theta(\bar x),
\]
contrary to $\mathcal B^+\models S_T$.

Consequently, there exists a model $\mathcal W\models T$
such that, upon identifying each $b\in B$ with the
interpretation of $c_b$,
\[
\mathcal B^+\prec\mathcal P_{\mathcal W}^{+}.
\]
Put
\[
\mathcal N=\Sk^{\mathcal W}(B)\prec\mathcal W.
\]
We claim that
\begin{equation}\label{eq:reconstructed-domain}
N\cap P^{\mathcal W}=B.
\end{equation}
The inclusion from right to left is immediate.
For the reverse inclusion, let
$d\in N\cap P^{\mathcal W}$. Choose an ambient
formula $\varphi(y,\bar x)$ and a tuple $\bar b$
from $B$ such that
\[
\mathcal W\models
\forall y\,(\varphi(y,\bar b)\leftrightarrow y=d).
\]
Since $d\in P^{\mathcal W}$,
\[
\mathcal P_{\mathcal W}^{+}
\models\exists!y\,R_\varphi(y,\bar b).
\]
By elementarity,
\[
\mathcal B^+\models\exists!y\,R_\varphi(y,\bar b).
\]
Let $e\in B$ be the witness. Then
\[
\mathcal W\models\varphi(e,\bar b),
\]
whence $e=d$, establishing
\eqref{eq:reconstructed-domain}.

It follows that
\[
P^{\mathcal N}=N\cap P^{\mathcal W}=B.
\]
Moreover, for each ambient formula $\psi(\bar x)$
and tuple $\bar b$ from $B$,
\begin{align*}
\mathcal P_{\mathcal N}^{+}
\models R_\psi(\bar b)
&\iff\mathcal N\models\psi(\bar b)\\
&\iff\mathcal W\models\psi(\bar b)\\
&\iff\mathcal P_{\mathcal W}^{+}
\models R_\psi(\bar b)\\
&\iff\mathcal B^+\models R_\psi(\bar b).
\end{align*}
This establishes (a). When the language includes
interpreted arithmetic operations, their graphs agree
by the same argument.

By elementarity of $\mathcal N$ in $\mathcal W$,
\[
\Sk^{\mathcal N}(B)
=\Sk^{\mathcal W}(B)=\mathcal N,
\]
establishing (b).

Finally, since the ambient language is countable,
counting formulas and finite tuples gives
\[
|B|\leq |N|
\leq\aleph_0\cdot|B|^{<\omega}
=\max(\aleph_0,|B|).
\]
This establishes (c). By the preceding term-model construction, $\mathcal K(\mathcal B^+)\cong\mathcal N$ via evaluation, which sends $\iota_B(b)$ to $b$. In particular, the designated domain of the reconstructed model is exactly $\iota_B[B]$, and $\iota_B$ is an isomorphism from $\mathcal B^+$ onto its induced expansion. Whenever literal equality of designated domains is desired, we transport the reconstructed model along a bijection extending $\iota_B^{-1}$.
\end{proof}

The central property of the reconstructed parent models
is that elementary embeddings between their induced expansions extend uniquely to elementary embeddings
between the parent models.

\begin{theorem}[Transfer]\label{transfer}
Let $\mathcal N_1,\mathcal N_2\models T$, with designated
domains $B_1,B_2$ and corresponding induced expansions
$\mathcal B_1^+,\mathcal B_2^+$. Suppose that $\mathcal N_i$
is definably generated by $B_i$ for $i=1,2$.
Then:
\begin{enumerate}[label=(\alph*)]
\item Every elementary embedding
\[
j:\mathcal B_1^+\longrightarrow\mathcal B_2^+
\]
extends uniquely to an elementary embedding
\[
\widehat j:\mathcal N_1\longrightarrow\mathcal N_2.
\]
Conversely, every elementary embedding between the
parent models restricts to an elementary embedding
between their induced expansions.

\item For every parent model $\mathcal N$ definably generated by its designated domain $B$ and induced expansion $\mathcal B^+$, the restriction map
\[
\begin{aligned}
\rho:\Aut(\mathcal N)
&\longrightarrow\Aut(\mathcal B^+),\\
h&\longmapsto h\restriction B
\end{aligned}
\]
is an isomorphism of topological groups.
\end{enumerate}
\end{theorem}

\begin{proof}
(a) Given $j$ and $a\in N_1$, choose a formula
$\varphi(y,\bar x)$ and a tuple $\bar b$ from $B_1$
such that $a$ is the unique solution of $\varphi(y,\bar b)$
in $\mathcal N_1$. Since
\[
\mathcal B_1^+\models R_{\exists!y\,\varphi(y,\bar x)}(\bar b),
\]
elementarity of $j$ ensures that $\varphi(y,j(\bar b))$
has a unique solution $a'$ in $\mathcal N_2$. Put $\widehat j(a)=a'$.

For any formula $\psi(y_1,\ldots,y_k)$ and elements $a_i$
uniquely defined by $\varphi_i(y_i,\bar b_i)$, put
\[
\chi(\bar x_1,\ldots,\bar x_k):=
\exists y_1\cdots\exists y_k
\left(\bigwedge_{i=1}^k\varphi_i(y_i,\bar x_i)
\land\psi(y_1,\ldots,y_k)\right).
\]
The truth of $\psi$ on the chosen elements is equivalent to
that of $R_\chi$ on their defining parameters in the induced expansion, both before and after applying $j$. Taking $\psi$
to be equality establishes independence of the defining formulas
and injectivity; the general case establishes elementarity.
Using the formula $y=x$ shows that $\widehat j\restriction B_1=j$.
Uniqueness follows from definable generation.

Conversely, suppose
$h:\mathcal N_1\longrightarrow\mathcal N_2$
is an elementary embedding. Since the designated
domains are defined by the same parameter-free
formula $P(x)$,
\[
h[B_1]\subseteq B_2.
\]
By \eqref{eq:kunen-translation},
\begin{align*}
\mathcal B_1^+\models\theta(\bar b)
&\iff\mathcal N_1\models\theta^\dagger(\bar b)\\
&\iff\mathcal N_2\models\theta^\dagger(h(\bar b))\\
&\iff\mathcal B_2^+\models\theta(h(\bar b)).
\end{align*}
Thus $h\restriction B_1$ is elementary.

(b) By (a), every automorphism $j$ of $\mathcal B^+$
extends uniquely to an elementary embedding
$\widehat j$ of $\mathcal N$ into itself. Applying
(a) to $j^{-1}$ and using uniqueness, we obtain
\[
\widehat{j^{-1}}\circ\widehat j
=\widehat j\circ\widehat{j^{-1}}
=\operatorname{id}_N.
\]
Hence $\widehat j\in\Aut(\mathcal N)$, and restriction
is a group isomorphism.

For each finite $F\subseteq B$,
\[
\rho^{-1}[\Aut(\mathcal B^+/F)]
=\Aut(\mathcal N/F).
\]
This establishes continuity of $\rho$ at the identity.

Conversely, given a finite $C\subseteq N$, choose
for each $c\in C$ a finite tuple $\bar b_c$ from $B$
from which $c$ is uniquely definable in $\mathcal N$.
Let $F\subseteq B$ be the finite set of all entries
of these tuples. Then
\[
\Aut(\mathcal N/F)\subseteq\Aut(\mathcal N/C),
\]
and consequently
\[
\rho^{-1}[\Aut(\mathcal B^+/F)]
\subseteq\Aut(\mathcal N/C).
\]
By the neighborhood bases described in
Subsection~\ref{subsec:topology}, the inverse of
$\rho$ is continuous at the identity. Since both
maps are homomorphisms, they are continuous
everywhere.
\end{proof}

The Transfer Theorem yields the following uniqueness result.

\begin{corollary}\label{kunen-categorical}
Any two parent models definably generated by their designated
domains and having the same induced expansion are isomorphic
via a unique isomorphism fixing the designated domain pointwise.
\end{corollary}

\begin{proof}
Apply Theorem~\ref{transfer} to the identity map on the common
induced expansion in both directions.
\end{proof}

\begin{remark}\label{kunen-functor}
In categorical terms, the Transfer Theorem shows that the induced expansion operation is full and faithful on the category of parent
models definably generated by their designated domains, with
 elementary embeddings as morphisms: every elementary embedding
between their induced expansions lifts uniquely to an elementary
embedding between the parent models. When $T$ has definable Skolem
functions, the reconstruction theorem ensures that every model of
$S_T$ arises, up to isomorphism, as such a induced expansion.
Thus in this case the induced expansion operation yields an
equivalence of categories.
\end{remark}

The following remark explores the extent to which
a parent model and its arithmetic induced expansion
are bi-interpretable.

\begin{remark}[Infinitary bi-interpretability]
\label{kunen-biinterpretability}
Suppose $\mathcal N$ interprets a model $\mathcal A$
of $\PA$ via a parameter-free, one-dimensional,
equality-preserving interpretation, and that
$\mathcal N$ is definably generated by the universe
of $\mathcal A$. Let $\mathcal A^+$ be the corresponding
induced expansion.
Then $\mathcal N$ and $\mathcal A^+$ are
bi-interpretable in $\mathcal L_{\omega_1,\omega}$.

The interpretation of $\mathcal A^+$ in $\mathcal N$
is first-order in the sense that the domain and each individual
symbol have first-order definitions; no uniform first-order
satisfaction predicate for the family of named relations is
asserted. For the reverse interpretation, enumerate the formulas
$\varphi_i(y,\bar x_i)$ of the parent language, with
$n_i=|\bar x_i|$, and let $i_0$ index the formula $y=x$.
Write $f_i^{\mathcal N}(\bar a)$ for the unique solution of
$\varphi_i(y,\bar a)$ whenever it exists uniquely.
Arithmetic coding represents the pair $\langle i,\bar a\rangle$
by an element of $A$, using the standard numeral for $i$.
Let $D_i(x)$ define the set of such codes for which
$R_{\exists!y\,\varphi_i(y,\bar x_i)}(\bar a)$ holds.
The interpretation domain is defined by
\[
D(x):=\bigvee_{i<\omega}D_i(x).
\]
For fixed $i,j<\omega$, let $Q_{ij}(x,y)$ express equality
of the unique solutions associated with the decoded tuples,
using the named predicate for
$\exists z(\varphi_i(z,\bar u)\land\varphi_j(z,\bar v))$.

The interpreted equality is therefore defined by
\[
E(x,y):=
\bigvee_{i,j<\omega}
\bigl(D_i(x)\land D_j(y)\land Q_{ij}(x,y)\bigr).
\]
Similarly, for each $k$-ary ambient relation symbol
$R$, let $Q^R_{\bar i}(\bar x)$ express
\[
R^{\mathcal N}
(f_{i_1}(\bar a_1),\ldots,f_{i_k}(\bar a_k)),
\]
on the corresponding coded tuples, using the named predicate for
$\exists y_1\cdots\exists y_k(\bigwedge_{\ell=1}^k
\varphi_{i_\ell}(y_\ell,\bar x_\ell)\land R(y_1,\ldots,y_k))$. The interpretation
of $R$ is defined by
\[
R^I(\bar x):=
\bigvee_{\bar i\in\omega^k}
\left(
\bigwedge_{\ell=1}^{k}D_{i_\ell}(x_\ell)
\land Q^R_{\bar i}(\bar x)
\right).
\]
Function symbols are treated through their graphs.

The resulting quotient structure is isomorphic to $\mathcal N$ via the evaluation map
\[
\begin{aligned}
\varepsilon:D/E&\longrightarrow N,\\
[\langle i,\bar a\rangle]_E
&\longmapsto f_i^{\mathcal N}(\bar a).
\end{aligned}
\]
Surjectivity follows from definable generation;
injectivity and preservation of the ambient relations
follow from the definitions of $E$ and $R^I$.

On the parent side, the graph of the evaluation map
on representatives is infinitarily definable by
\[
\bigvee_{i<\omega}
\exists\bar a\,
\bigl(x=\langle i,\bar a\rangle
\land\varphi_i(y,\bar a)\land\exists!z\,\varphi_i(z,\bar a)\bigr),
\]
where the coding and decoding are interpreted in the
parent model, and the quantified tuples range over
the designated domain.

On the child side, the isomorphism witnessing the
other composite is
\[
a\longmapsto[\langle i_0,a\rangle]_E.
\]
Indeed, under $\varepsilon$,
\[
[\langle i_0,a\rangle]_E
\longmapsto f_{i_0}^{\mathcal N}(a)=a.
\]
Thus the composite recovers the designated arithmetic
domain, its arithmetic operations, and all its named
predicates $R_\varphi$. The two composite isomorphisms
establish the asserted infinitary bi-interpretability.

By a result of Visser, no model of $\PA$ is
bi-interpretable with a model of $\ZF$; indeed, neither
can be a parametric retract of the other
(see \cite[Corollaries~9.6 and~9.8]{Visser2006} and
\cite[proof of Theorem~18]{EL}).
A related argument shows that the passage to
the induced expansion does not overcome the first-order
obstruction. More precisely, suppose
$\mathcal N\models\ZF+\mathsf{V=OD}$ and $\mathcal A$
is its canonically interpreted arithmetic structure.
Suppose that there were a parameter-free first-order
interpretation of $\mathcal N$ in $\mathcal A^+$. Composing it with the canonical
interpretation of $\mathcal A^+$ in $\mathcal N$
would yield a parameter-free definable set-sized
structure $\mathcal D$ within $\mathcal N$,
externally isomorphic to $\mathcal N$.

Using the definability of satisfaction for set-sized
structures, the predicate
\[
\operatorname{Tr}(n)
\quad\Longleftrightarrow\quad
\operatorname{Sat}(\mathcal D,n,\varnothing)
\]
would be parameter-free definable in $\mathcal N$
and satisfy, for every standard sentence $\sigma$
of set theory,
\[
\mathcal N\models
\operatorname{Tr}(\ulcorner\sigma\urcorner)
\quad\Longleftrightarrow\quad
\mathcal N\models\sigma,
\]
contradicting Tarski's undefinability theorem.
\end{remark}

\section{Kunen expansions meet Schmerl realizations}\label{sec:realization}
Theorems~\ref{kunen-reconstruction} and~\ref{transfer} reduce the realization problem for accommodating theories with definable Skolem functions to realization in expanded arithmetic. Given a model $\mathcal A^+\models S_T$ supplied by arithmetic realization, we obtain
\[
\mathcal K(\mathcal A^+)\models T,
\qquad
\Aut(\mathcal K(\mathcal A^+))\cong_{\mathrm{top}}\Aut(\mathcal A^+).
\]
The reconstruction theorem supplies definably generated parent models, to which the Transfer Theorem applies without any further Skolem hypothesis. We begin with the arithmetic results needed for these applications.

Schmerl's characterization of the spectrum of $\PA$ (Theorem~\ref{schmerl-spectrum}) extends to expanded languages. For Theorem A(2) we need the following refinement, which permits us to prescribe the automorphism group of an elementary end extension of a given expanded model of arithmetic.
\begin{theorem}[{Kossak--Schmerl \cite[Theorem 5.4.4, p.~150]{KS}}]\label{ks}
Let $G$ be right-orderable and let $\mathcal A_0^{+}$ be a model of $\PA^{*}$. There exists an elementary end extension $\mathcal A^{+}$ of $\mathcal A_0^{+}$ such that
\[
\Aut(\mathcal A^{+})\cong G.
\]
\end{theorem}
To obtain realization results that also control the topology of the automorphism group, we use the following topological form of Schmerl's Theorem A, which follows from the construction in the proof of Theorem~5.4.4 of Kossak and Schmerl.\footnote{Schmerl states his theorem for finite expanded languages (see \cite[p.~1251]{Schmerl}). The countable-language version follows from the proof of \cite[Theorem~5.4.4, pp.~150--151]{KS}, using the coding in the proof of \cite[Theorem~5.4.3, pp.~149--150]{KS} and starting with a countable model of $T$. The isomorphism in \cite[Proposition~5.1.2]{KS} is topological: restriction to a solid basis is continuous, and continuity of its inverse follows from finite Skolem dependence. The same argument applies in the proof of Theorem~5.4.4, since the auxiliary base and the additional parameter $c$ are fixed pointwise by all automorphisms. Restriction from the coded ordered structure to the original ordered structure is also a topological isomorphism, its inverse being extension by the identity on the auxiliary well-ordered block.}

\begin{theorem}[{Schmerl; see \cite[Theorem A]{Schmerl} and \cite[proof of Theorem~5.4.4]{KS}}]\label{main}
Let $T$ be a consistent extension of $\PA^{*}$ in the same language. For every nonempty linearly ordered structure $\mathcal K$ in a countable language, there exists $\mathcal A^{+}\models T$ such that
\[
\Aut(\mathcal K)\cong_{\mathrm{top}}\Aut(\mathcal A^{+}),
\qquad |A|=\max(\aleph_0,|K|).
\]
\end{theorem}

\begin{corollary}\label{equal}
For every consistent extension $T$ of $\PA^*$ in the same language, the automorphism-group spectrum of $T$ consists precisely of the right-orderable groups. Thus
\[
\G(T)=\G(\PA)=\Grpr.
\]
\end{corollary}
\begin{proof}
For any $\mathcal A^+\models T$, restriction to the arithmetic reduct embeds $\Aut(\mathcal A^+)$ into $\Aut(\mathcal A)$, which is right-orderable by Theorem~\ref{orderability}. Thus $\G(T)\subseteq\Grpr$. Conversely, given $G\in\Grpr$, choose $\mathcal A_0^+\models T$ and apply Theorem~\ref{ks} to obtain $\mathcal A_0^+\prec_{\mathrm{end}}\mathcal A^+$ with $\Aut(\mathcal A^+)\cong G$. Since $\mathcal A^+\models T$ by elementarity, $G\in\G(T)$.
\end{proof}

\subsection{General realization theorems}
For the realization theorems below we use parameter-free, one-dimensional, equality-preserving interpretations witnessing accommodation. Accommodation allows us to apply the arithmetic realization theorems to Kunen expansions, while definable Skolem functions allow us to reconstruct parent models without changing the automorphism groups. A definable linear ordering of the parent model serves a different purpose: it ensures that every automorphism group of a model of the ambient theory is right-orderable.

\begin{theorem}\label{general}
Let $T$ be a theory in a countable language with definable Skolem functions accommodating $\PA$ via a parameter-free, one-dimensional, equality-preserving interpretation. Then
\[
\Grpr\subseteq\G(T).
\]
If $T$ also has a parameter-free definable linear ordering of its entire domain, then equality holds.
\end{theorem}
\begin{proof}
Fix a parent model $\mathcal M\models T$. By Proposition~\ref{accommodation-induction}, its arithmetic Kunen expansion satisfies $\PA^*$; consequently $S_T\supseteq\PA^*$. Given $G\in\Grpr$, Corollary~\ref{equal} yields $\mathcal A^+\models S_T$ with $\Aut(\mathcal A^+)\cong G$. Put $\mathcal N=\mathcal K(\mathcal A^+)$, the associated Kunen model. By Theorem~\ref{kunen-reconstruction}, $\mathcal N\models T$ is generated under definable Skolem functions by its arithmetic domain, and $\mathcal A_{\mathcal N}^+\cong\mathcal A^+$. Theorem~\ref{transfer} gives
\[
\Aut(\mathcal N)\cong\Aut(\mathcal A^+)\cong G.
\]

If $T$ has a parameter-free definable linear ordering $\triangleleft$ of its entire domain, then for every $\mathcal M\models T$,
\[
\Aut(\mathcal M)\leq\Aut(M,\triangleleft^{\mathcal M}).
\]
By Theorem~\ref{orderability}, $\G(T)\subseteq\Grpr$.
\end{proof}

\begin{remark}\label{rem:skolem-choice}
Accommodation alone does not ensure faithful restriction to the arithmetic domain. Cohen's theorem \cite[Theorem 1]{Cohen} gives a model $\mathcal M\models\ZF$ with a nonidentity ordinal-fixing involution. Since $\ZF$ accommodates $\PA$ via the canonical interpretation on $\omega$, this involution lies in the kernel of restriction to $\mathcal A_{\mathcal M}^+$. By Proposition~\ref{faithful}, Choice fails in $\mathcal M$.

Indeed, definable Skolem functions imply Choice over $\ZF$, even if parameters are allowed in their definitions. To see this, let $f$ be a definable Skolem function for the formula $x\in y$. Working within a model of $\ZF$ in which such a function is definable, let $a$ be any set of nonempty sets. Replacement yields the set
\[
\{\langle y,f(y)\rangle:y\in a\},
\]
which is a choice function on $a$, since $f(y)\in y$ for every $y\in a$. Thus Cohen's model cannot have definable Skolem functions.
\end{remark}

The topological form of transfer allows us to retain the topology of the automorphism group of an ordered structure, as well as to control the cardinality of the resulting parent model.
\begin{theorem}\label{general-top}
Suppose $T$ satisfies the hypotheses of Theorem~\ref{general} concerning accommodation and definable Skolem functions. For every nonempty linearly ordered structure $\mathcal K$ in a countable language, there exists $\mathcal N\models T$ such that
\[
\Aut(\mathcal N)\cong_{\mathrm{top}}\Aut(\mathcal K),\qquad |N|=\max(\aleph_0,|K|).
\]
\end{theorem}
\begin{proof}
As in the proof of Theorem~\ref{general}, $S_T\supseteq\PA^*$. Apply Theorem~\ref{main} to $S_T$ and $\mathcal K$ to obtain $\mathcal A^+\models S_T$ with
\[
\Aut(\mathcal A^+)\cong_{\mathrm{top}}\Aut(\mathcal K),\qquad |A|=\max(\aleph_0,|K|).
\]
Put $\mathcal N=\mathcal K(\mathcal A^+)$, the associated Kunen model. By Theorem~\ref{kunen-reconstruction}, this is a model of $T$ definably generated by its arithmetic domain with $\mathcal A_{\mathcal N}^+\cong\mathcal A^+$ and $|N|=|A|$. By Theorem~\ref{transfer}, $\Aut(\mathcal N)\cong_{\mathrm{top}}\Aut(\mathcal A^+)$.
\end{proof}
\subsection{Applications to set theory}
\label{sec:applications}
The canonical interpretation of arithmetic on $\omega$ witnesses accommodation over $\ZF$ by Separation. In models of $\ZF+\mathsf{V=OD}$, the canonical ordering by least ordinal-definition codes is a parameter-free definable set-like global well-order. Choosing least witnesses therefore gives definable Skolem functions. Moreover, over $\ZF$, $\mathsf{V=OD}$ implies both $\mathsf{AC}$ and $\mathsf{LM}$: the global well-order yields Choice, and the ordinal-definition codes give a parameter-free definable injection into the class of sets of ordinals, which is equivalent to $\mathsf{LM}$ by \cite[Theorem~2.1]{LM}. These observations allow us to apply both general realization theorems.

\begin{corollary}\label{od-spectrum}
For every consistent extension $T$ of $\ZF+\mathsf{V=OD}$ in $\Lset$,
\[
\G(T)=\Grpr=\G(\PA).
\]
In particular, the same holds for every consistent extension of $\ZF+\mathsf{V=L}$.
\end{corollary}
\begin{proof}
Apply Theorem~\ref{general}, using the definable global well-order described above. The last assertion follows since $\mathsf{V=L}$ implies $\mathsf{V=OD}$.
\end{proof}

\begin{corollary}\label{od-top}
Let $T$ be a consistent extension of $\ZF+\mathsf{V=OD}$ in $\Lset$. For every nonempty linearly ordered structure $\mathcal K$ in a countable language, there is $\mathcal N\models T$ such that
\[
\Aut(\mathcal N)\cong_{\mathrm{top}}\Aut(\mathcal K),\qquad |N|=\max(\aleph_0,|K|).
\]
\end{corollary}
\begin{proof}
Apply Theorem~\ref{general-top} under the same hypotheses.
\end{proof}

A definable global well-order is not needed to ensure right-orderability: a faithful action on the ordinals suffices.\footnote{The assertion of Proposition~\ref{faithful} for \textup{\textsf{ZFC}} is due to Friedman, as reported in \cite[p.~326]{Cohen}. The extension to \textup{\textsf{LM}} is implicit in the proof of \cite[Theorem 3.6(2)]{Paris}, where it is used to establish rigidity of Paris models (models whose ordinals are all parameter-free definable).}

\begin{proposition}\label{faithful}
If $\mathcal M\models\mathsf{ZFC}$ or $\mathcal M\models\ZF+\mathsf{LM}$, restriction induces an injective homomorphism
\[
\Aut(\mathcal M)\longrightarrow\Aut(\operatorname{Ord}^{\mathcal M},\in^{\mathcal M}).
\]
\end{proposition}
\begin{proof}
Let $j\in\Aut(\mathcal M)$ fix every ordinal of $\mathcal M$. By Extensionality, $j$ fixes every set of ordinals of $\mathcal M$.

Suppose first that $\mathcal M\models\mathsf{ZFC}$. Working within $\mathcal M$, for any $a\in M$ choose a bijection from an ordinal $\alpha$ onto the transitive closure of $\{a\}$. Transport $\in$ to a well-founded extensional relation on $\alpha$, and let $\gamma<\alpha$ correspond to $a$. The relation is coded by a set of ordinals $s$ using the canonical ordinal pairing function. By the Mostowski collapse theorem, $a$ is uniquely definable from $s$ and $\gamma$. Since both parameters are fixed by $j$, so is $a$. The collapse is internal to $\mathcal M$; no external well-foundedness is required.

For models satisfying the Leibniz--Mycielski axiom, we use an equivalent formulation. By \cite[Theorem~2.1]{LM}, see also \cite[Theorem 3.4]{Paris}, \textup{\textsf{LM}} is equivalent over \textup{\textsf{ZF}} to the existence of a parameter-free definable injection $F$ from the universe into the class of sets of ordinals. Hence for every $a\in M$,
\[
F(j(a))=j(F(a))=F(a),
\]
so $j(a)=a$. In either case $j$ is the identity.
\end{proof}

\begin{corollary}\label{choice-spectrum}
We have
\[
\G(\mathsf{ZFC})=\G(\ZF+\mathsf{LM})=\Grpr\subsetneq\G(\ZF).
\]
\end{corollary}
\begin{proof}
Proposition~\ref{faithful} and Theorem~\ref{orderability} give the inclusions in $\Grpr$, since the ordinals of a model of \textup{\textsf{ZF}} are linearly ordered. The reverse inclusions follow from Corollary~\ref{od-spectrum} applied to $\ZF+\mathsf{V=L}$, which is consistent relative to $\ZF$ by constructibility and implies both \textup{\textsf{AC}} and \textup{\textsf{LM}}.

For strictness, use Cohen's theorem \cite[Theorem 1]{Cohen}: every countable model of \textup{\textsf{ZF}} with an externally ill-founded ordinal gives rise to a model of \textup{\textsf{ZF}} with the same ordinals and a nonidentity automorphism of order two. See also \cite[proof of Theorem 3.6(3)]{Paris} and the alternative proof in \cite{An1} (see also \cite{HodgesReview}). In particular, by compactness and downward L\"owenheim--Skolem, there is a countable model of \textup{\textsf{ZF}} with nonstandard natural numbers to which this result applies. The full automorphism group of the resulting model has torsion, and therefore is not right-orderable.
\end{proof}

Corollaries~\ref{od-spectrum} and~\ref{choice-spectrum} establish Theorem A(1). 

For an arbitrary consistent completion $T$ of \textup{\textsf{ZFC}} or of $\ZF+\mathsf{LM}$, Proposition~\ref{faithful} and Theorem~\ref{orderability} give $\G(T)\subseteq\Grpr$. Corollary~\ref{od-spectrum} gives equality when $T$ extends $\ZF+\mathsf{V=OD}$; the argument does not establish equality for arbitrary completions. Nor does it characterize the entire spectrum of $\ZF$.

\subsection{Elementary end embeddings}
\label{sec:endsec}
The Kossak--Schmerl theorem (Theorem~\ref{ks}) allows us to prescribe the automorphism group of a model while retaining the arithmetic Kunen expansion of a given parent model as an elementary initial segment. Reconstruction then yields an elementary embedding of the Skolem hull of the original arithmetic domain, and Ackermann coding yields an elementary end embedding of the corresponding hereditarily finite structures.

\begin{theorem}\label{hfend}
Suppose $\mathcal M\models\ZF+\mathsf{V=OD}$ and $G$ is a right-orderable group. There is a model $\mathcal N\equiv\mathcal M$ such that:
\begin{enumerate}[label=(\alph*)]
\item $\Aut(\mathcal N)\cong G$.
\medskip
\item There is an elementary end embedding
\[
j:\mathcal A_{\mathcal M}^{+}\longrightarrow\mathcal A_{\mathcal N}^{+}.
\]
\item $\mathcal N=\Sk^{\mathcal N}(\omega^{\mathcal N})$, and $j$ extends uniquely to an elementary embedding
\[
\widehat j:\Sk^{\mathcal M}(\omega^{\mathcal M})\longrightarrow\mathcal N.
\]
\item Via the canonical Ackermann bijections, $j$ induces an elementary end embedding
\[
(V_\omega^{\mathcal M},\in^{\mathcal M})\longrightarrow
(V_\omega^{\mathcal N},\in^{\mathcal N}).
\]
\end{enumerate}
\end{theorem}
\begin{proof}
Put $T=\operatorname{Th}(\mathcal M)$ and take the canonical interpretation of arithmetic on $\omega$ as the designated interpretation. Since $\ZF$ accommodates $\PA$ via this interpretation, Proposition~\ref{accommodation-induction} gives $\mathcal A_{\mathcal M}^{+}\models\PA^*$. By Theorem~\ref{ks}, there is an elementary end extension
\[
\mathcal A_{\mathcal M}^{+}\prec_{\mathrm{end}}\mathcal A^{+},
\qquad \Aut(\mathcal A^{+})\cong G.
\]
By elementarity, $\mathcal A^{+}\models S_T$. Since $T$ has definable Skolem functions, take $\mathcal N$ to be the Kunen model $\mathcal K(\mathcal A^+)$ transported along the canonical identification of its designated domain with $A$. By Theorem~\ref{kunen-reconstruction}, $\mathcal N\models T$ and
\[
\mathcal A_{\mathcal N}^{+}=\mathcal A^{+},
\qquad \mathcal N=\Sk^{\mathcal N}(\omega^{\mathcal N}).
\]
Theorem~\ref{transfer} yields
\[
\Aut(\mathcal N)\cong\Aut(\mathcal A^{+})\cong G.
\]
Let $j$ be the inclusion of $\mathcal A_{\mathcal M}^{+}$ into $\mathcal A_{\mathcal N}^{+}$. This establishes (a) and (b).

Put $\mathcal M^*=\Sk^{\mathcal M}(\omega^{\mathcal M})$. By Proposition~\ref{kunen-internal-reconstruction},
\[
\omega^{\mathcal M^*}=\omega^{\mathcal M},
\qquad \mathcal A_{\mathcal M^*}^{+}=\mathcal A_{\mathcal M}^{+}.
\]
Both $\mathcal M^*$ and $\mathcal N$ are definably generated by their arithmetic domains. By Theorem~\ref{transfer}, $j$ extends uniquely to the elementary embedding $\widehat j$ in (c).

For (d), let $c_{\mathcal K}:V_\omega^{\mathcal K}\to\omega^{\mathcal K}$ denote the canonical Ackermann bijection for $\mathcal K\in\{\mathcal M,\mathcal N\}$. Under this coding, $\in^{\mathcal K}$ corresponds to the arithmetically definable relation
\[
u\in_{\mathrm{Ack}}v
\quad\Longleftrightarrow\quad
\text{the $u$th binary digit of $v$ is $1$}.
\]
Define
\[
j_{\mathrm{HF}}=c_{\mathcal N}^{-1}\circ j\circ c_{\mathcal M}.
\]
The arithmetic reduct of $j$ is elementary. Thus, for each $\Lset$-formula $\varphi$, with arithmetic translation $\varphi^{\mathrm{Ack}}$, and each tuple $\bar a$ from $V_\omega^{\mathcal M}$,
\begin{align*}
(V_\omega^{\mathcal M},\in^{\mathcal M})\models\varphi(\bar a)
&\iff\mathcal A_{\mathcal M}\models\varphi^{\mathrm{Ack}}(c_{\mathcal M}(\bar a))\\
&\iff\mathcal A_{\mathcal N}\models\varphi^{\mathrm{Ack}}(j(c_{\mathcal M}(\bar a)))\\
&\iff(V_\omega^{\mathcal N},\in^{\mathcal N})\models\varphi(j_{\mathrm{HF}}(\bar a)).
\end{align*}
Hence $j_{\mathrm{HF}}$ is elementary.

To establish the end-embedding property, let $a\in V_\omega^{\mathcal M}$ and suppose $d\in^{\mathcal N}j_{\mathrm{HF}}(a)$. Then $d\in V_\omega^{\mathcal N}$ and
\[
c_{\mathcal N}(d)\in_{\mathrm{Ack}}j(c_{\mathcal M}(a)).
\]
Since $\PA$ proves $u\in_{\mathrm{Ack}}v\Rightarrow u<v$, we have
\[
c_{\mathcal N}(d)<^{\mathcal A_{\mathcal N}}j(c_{\mathcal M}(a)).
\]
The image of $j$ is an initial segment, so there is $b\in\omega^{\mathcal M}$ with $c_{\mathcal N}(d)=j(b)$. Consequently,
\[
d=j_{\mathrm{HF}}(c_{\mathcal M}^{-1}(b)).
\]
By elementarity, $c_{\mathcal M}^{-1}(b)\in^{\mathcal M}a$. This establishes (d).
\end{proof}

\section{Tzouvaras meets Enayat--Schmerl--Visser}
\label{sec:nofoundation}
We now turn to Theorem B, whose proof does not make
use of the machinery of Sections~\ref{sec:arithmetic-transfer}
and~\ref{sec:realization}. Our construction extends the
method underlying Theorem~\ref{esv-spectrum} and has
antecedents in Tzouvaras \cite[Section~4]{Tz} and
Enayat, Schmerl and Visser \cite{ESV}. In both constructions
one successively adjoins representatives of missing subsets:
in \cite{ESV} the subsets are finite and the iteration is
through $\omega$, whereas in \cite[Section~4]{Tz} arbitrary
subsets are used and the iteration is transfinite.
The transfinite construction was also described by
the author in \cite{MO}, under the name von Neumann
completion. Our variant differs from Tzouvaras's in
that new vertices are specified directly by tagging,
rather than allocated by means of a global well-order.
As shown in Appendix~\ref{sec:class-completions}, this
eliminates the need for Global Choice even in the
proper-class setting. In the appendix we also revisit
the question concerning descending membership chains
in models with well-founded ordinals that motivated
the construction in \cite{MO}.

Tzouvaras completion also provides a versatile alternative
to the Rieger--Bernays permutation method (see, e.g.,
\cite[Chapter~3]{Forster}) for establishing many classical
results concerning set theory without Foundation. For instance,
one can prescribe finite membership cycles or sets of Quine
atoms by choosing suitable extensional seed digraphs.

We combine the completion construction developed in Subsection~\ref{subsec:completion} with the coding digraph of Enayat, Schmerl and Visser \cite[Theorem 4.2]{ESV}, described in Subsection~\ref{subsec:recovery}.

\begin{theorem}\label{inacc}
Suppose $\kappa$ is a strongly inaccessible cardinal. For every group $H$ of cardinality less than $\kappa$, there is a structure $\mathcal N=(N,E)$ such that
\[
\mathcal N\models\mathsf{ZFC}^{-\mathsf{Found}},\qquad |N|=\kappa,
\qquad \Aut(\mathcal N)\cong H.
\]
The externally well-founded part of $\mathcal N$ is isomorphic to $(V_\kappa,\in)$.
\end{theorem}

\subsection{Tzouvaras completions}\label{subsec:completion}
To obtain a model of set theory from an extensional digraph, we must ensure that sufficiently many subsets of its universe are represented by fibres. We first develop the construction and its automorphism extension property without any cardinal hypothesis. Let $\mathcal G=(G,E)$ be extensional. Extensionality says that $x\mapsto x_E$ is injective; by Cantor's theorem it cannot be surjective. 

\begin{definition}\label{def:completion}
Let $\mathcal G=(G,E)$ be an extensional digraph.
\begin{enumerate}[label=(\alph*)]
\item
We call the members of
\[
D(\mathcal G):=\mathcal P(G)\setminus\{x_E:x\in G\}
\]
its \emph{deficiency sets}. We adjoin a vertex for each such subset, with precisely that subset as its fibre, and iterate the operation transfinitely.

\item To distinguish original vertices from newly adjoined ones, replace $\mathcal G$ by an isomorphic copy whose vertices have the form $\langle0,x\rangle$. We can then use $\langle1,S\rangle$ to represent a deficiency set $S$. We call this device \emph{tagging}. 
\item The \emph{Tzouvaras completion} is obtained by putting $\mathcal G_0=\mathcal G$ and, at successor stages,
\[
\begin{aligned}
D_\alpha&=D(\mathcal G_\alpha),\\
G_{\alpha+1}&=G_\alpha\cup(\{1\}\times D_\alpha),\\
E_{\alpha+1}&=E_\alpha\cup\{\langle x,\langle1,S\rangle\rangle:S\in D_\alpha,\ x\in S\}.
\end{aligned}
\]
At each nonzero limit ordinal $\beta$, put
\[
G_\beta=\bigcup_{\alpha<\beta}G_\alpha,\qquad
E_\beta=\bigcup_{\alpha<\beta}E_\alpha.
\]
We write $\mathcal G_\alpha=(G_\alpha,E_\alpha)$ at every stage.
\item For an infinite cardinal $\kappa$, the \emph{$\kappa$-bounded version} of the construction uses only deficiency sets of cardinality less than $\kappa$. We use the same stage notation for this version when the bound is specified.
\end{enumerate}
\end{definition}
By transfinite induction the digraphs are extensional and each end extends its predecessors. Indeed, a proposed new vertex $\langle1,S\rangle$ is not original by tagging; if introduced at an earlier stage, its fibre would still be $S$ since subsequent stages leave existing fibres unchanged, contradicting $S\in D_\alpha$. Thus the vertices are genuinely new, and
\[
 x_{E_{\alpha+1}}=x_{E_\alpha}\ (x\in G_\alpha),\qquad
 (\langle1,S\rangle)_{E_{\alpha+1}}=S\ (S\in D_\alpha).
\]
For $x\in G_\alpha$ and $S\in D_\alpha$, we have
\[
x_{E_{\alpha+1}}=x_{E_\alpha}\ne S=(\langle1,S\rangle)_{E_{\alpha+1}}.
\]
Distinct deficiency sets give distinct new fibres. These equalities show that each successor stage is extensional and end extends its predecessor; the same properties hold at limits. In particular $D_\alpha\cap D_\beta=\varnothing$ for $\alpha<\beta$. Write
\[
\mathcal G_{<\kappa}=\bigcup_{\alpha<\kappa}\mathcal G_\alpha,\qquad G_{<\kappa}=\bigcup_{\alpha<\kappa}G_\alpha.
\]

\begin{proposition}\label{completion-aut}
For every ordinal $\lambda$, each automorphism of $\mathcal G$ extends uniquely to an automorphism of $\mathcal G_\lambda$. Restriction identifies the setwise stabilizer of $G$ in $\Aut(\mathcal G_\lambda)$ with $\Aut(\mathcal G)$ as abstract groups.
\end{proposition}
\begin{proof}
Starting with $h_0=h\in\Aut(\mathcal G)$, define
\[
h_{\alpha+1}(x)=\begin{cases}
h_\alpha(x)&x\in G_\alpha,\\
\langle1,h_\alpha[S]\rangle&x=\langle1,S\rangle,\ S\in D_\alpha.
\end{cases}
\]
At nonzero limit stages put $h_\beta=\bigcup_{\alpha<\beta}h_\alpha$. If $h_\alpha\in\Aut(\mathcal G_\alpha)$, then
\[
 h_\alpha[x_{E_\alpha}]=(h_\alpha(x))_{E_\alpha}
 \qquad(x\in G_\alpha).
\]
Consequently $S\in D_\alpha\iff h_\alpha[S]\in D_\alpha$ for $S\subseteq G_\alpha$, so the successor map is bijective. For $x\in G_\alpha$ and $S\in D_\alpha$,
\[
xE_{\alpha+1}\langle1,S\rangle\iff x\in S\iff h_\alpha(x)\in h_\alpha[S]\iff h_\alpha(x)E_{\alpha+1}\langle1,h_\alpha[S]\rangle.
\]
Old edges are preserved and no new vertex belongs to an old fibre or another new fibre. Thus the successor map is an automorphism, and the limit assertion follows by unions. For uniqueness, let $k$ be an automorphism of $\mathcal G_\lambda$ extending $h$. If $k\restriction G_\alpha=h_\alpha$ and $\alpha<\lambda$, then, writing $E'=E_\lambda$, for every $S\in D_\alpha$,
\[
(k(\langle1,S\rangle))_{E'}=k[S]
=h_\alpha[S]=(\langle1,h_\alpha[S]\rangle)_{E'}.
\]
Extensionality of $\mathcal G_\lambda$ therefore gives $k\restriction G_{\alpha+1}=h_{\alpha+1}$. Transfinite induction proves uniqueness even among automorphisms of the final structure. Restriction and extension are inverse homomorphisms on the indicated stabilizer.
\end{proof}

Definition~\ref{def:completion} and Proposition~\ref{completion-aut} apply at every ordinal stage. To obtain a set-sized model of set theory, we now impose a cardinal hypothesis that ensures that all sufficiently small subsets are represented, while every fibre remains small.

Proposition~\ref{completion-aut} applies verbatim to the bounded version of Definition~\ref{def:completion}, since automorphisms preserve cardinalities.

\begin{lemma}\label{bounded}
Suppose $\kappa$ is strongly inaccessible, $|G|\leq\kappa$, and $|x_E|<\kappa$ for every $x\in G$. The $\kappa$-bounded completion through $\kappa$ stages is extensional, has cardinality $\kappa$, and its fibres are precisely the subsets of its universe of cardinality less than $\kappa$, each represented exactly once.
\end{lemma}
\begin{proof}
Since $\kappa^{<\kappa}=\kappa$ by regularity and the strong limit property, induction gives $|G_\alpha|\leq\kappa$ for $\alpha\leq\kappa$. All fibres remain small by end extension. Every subset of the final universe of size less than $\kappa$ is contained in a bounded stage, by regularity, and is represented at the next stage if not already represented. If the final universe had size less than $\kappa$, the fibre map would be a bijection onto its entire power set, contradicting Cantor's theorem. Thus its size is $\kappa$. When $|G|<\kappa$, all stages below $\kappa$ have size less than $\kappa$ by inaccessibility, and the bounded and unrestricted constructions agree through $\kappa$.
\end{proof}

Rieger's theorem on full extensional systems gives sufficient conditions for a digraph representing all ambient set-sized subsets of its universe to satisfy set theory without Foundation. Lemma~\ref{axioms} establishes a bounded counterpart of this theorem, based on the presentation by Aczel \cite[Appendix B]{Ac}.
\begin{lemma}\label{axioms}
Under the hypotheses of Lemma~\ref{bounded}, $\mathcal G_{<\kappa}\models\mathsf{ZFC}\setminus\{\mathsf{Found}\}$.
\end{lemma}
\begin{proof}
Put $\mathcal N=\mathcal G_{<\kappa}$. For $S\subseteq N$ of size less than $\kappa$, denote its unique representative by $[S]$. Extensionality, Empty Set, Pairing and Separation follow directly from Lemma~\ref{bounded}. For Infinity, recursively put $u_0=[\varnothing]$ and $u_{n+1}=[\{u_0,\ldots,u_n\}]$; the set $U=\{u_n:n<\omega\}$ has cardinality less than $\kappa$. Since $(u_n)_E=\{u_i:i<n\}$, the internal successor of $u_n$ is $u_{n+1}$. Thus $[U]$ contains the empty set and is closed under internal successor, verifying Infinity.

For Union, given $a\in N$, regularity ensures that $S=\bigcup_{b\in a_E}b_E$ has size less than $\kappa$, and
\[
xE[S]\iff\exists b\,(xEb\land bEa).
\]
Thus $[S]$ represents the union of $a$. For Power Set, the collection $S=\{[X]:X\subseteq a_E\}$ has cardinality $2^{|a_E|}<\kappa$ by the strong limit property. Since $bE[S]\iff b_E\subseteq a_E$, its representative is the power set of $a$ in $\mathcal N$.

For Replacement, if a formula defines a function on $a_E$ in $\mathcal N$, its external range has cardinality at most $|a_E|<\kappa$. Writing $S=\{y\in N:\exists x\in a_E\ \mathcal N\models\varphi(x,y,\bar p)\}$, we have
\[
\mathcal N\models\forall y\,(y\in[S]\leftrightarrow\exists x\in a\,\varphi(x,y,\bar p)).
\]
This verifies Replacement. For Collection, suppose $\mathcal N\models\forall x\in a\,\exists y\,\varphi(x,y,\bar p)$. By ambient Choice select, for each $x\in a_E$, a witness $y_x\in N$ with $\mathcal N\models\varphi(x,y_x,\bar p)$. Put $C=\{y_x:x\in a_E\}$. Since $|C|\leq |a_E|<\kappa$, Lemma~\ref{bounded} gives its representative $[C]$, and
\[
\mathcal N\models\forall x\in a\,\exists y\in[C]\,\varphi(x,y,\bar p).
\]
This verifies Collection.

For Choice, choose an external well-ordering $R$ of $a_E$. Using internal Kuratowski ordered pairs, form the set $S_R=\{\langle x,y\rangle^{\mathcal N}:xRy\}$. Its cardinality is at most $|a_E|^2<\kappa$, so $[S_R]$ represents a relation on $a$ in $\mathcal N$. Every internally represented nonempty subset of $a$ corresponds to a nonempty subset of $a_E$ and therefore has an $R$-least element. Thus $\mathcal N$ regards this relation as a well-ordering of $a$, establishing Choice.
\end{proof}

\begin{proposition}\label{wf-completion}
Under the hypotheses of Lemma~\ref{bounded}, the externally well-founded part of the completion is isomorphic to $(V_\kappa,\in)$. In particular, the completion is ordinal-standard of ordinal height $\kappa$.
\end{proposition}
\begin{proof}
The externally well-founded part is downward closed and extensional, and hence admits a Mostowski collapse. By well-founded induction and regularity, every collapsed element has rank less than $\kappa$, since its fibre has size less than $\kappa$. Conversely, induction on rank shows that every element of $V_\kappa$ is in the range of the collapse: its members have already been represented, and their representatives form a set of size less than $\kappa$ by inaccessibility. Its unique representative is externally well-founded. Every internal ordinal is externally well-founded, since an external descending sequence in it would have a represented range contradicting its internal well-ordering. Conversely, the inverse image of each ordinal below $\kappa$ under the collapse is an internal ordinal, since the well-founded part is downward closed and the ambient ordinal ordering is a well-order. Thus the internal ordinals correspond exactly to the ordinals below $\kappa$.
\end{proof}

\begin{remark}
Conversely, if an infinite $\Lset$-structure of cardinality $\kappa$ has its fibres in bijective correspondence with all subsets of its universe of size less than $\kappa$, and satisfies Power Set, Union and Infinity, then $\kappa$ is strongly inaccessible. Power Set gives $2^\lambda<\kappa$ for every cardinal $\lambda<\kappa$. If $\kappa$ were singular, a family of subsets of sizes cofinal in $\kappa$, indexed by $\operatorname{cf}(\kappa)$, would have its representatives collected in a single fibre, whose union would have size $\kappa$, a contradiction. Infinity excludes $\aleph_0$.
\end{remark}

\subsection{From graphs to sets and back}\label{subsec:recovery}
Proposition~\ref{completion-aut} ensures that automorphisms of an extensional digraph extend uniquely to its completion. To ensure that these are all the automorphisms of the completion, we construct a digraph that can be recovered there without parameters.
Fix a strongly inaccessible cardinal $\kappa$ and a group $H$ of cardinality less than $\kappa$. By Theorem~\ref{graph-representation}, choose a simple undirected graph $\mathcal A=(A,F)$ with $\Aut(\mathcal A)\cong H$ and $|A|<\kappa$. We may take $A$ nonempty. The coding device of Enayat, Schmerl and Visser replaces each vertex by an infinite descending $E$-sequence, whose positions are distinguished by finite ordinals. A separate descending chain for each edge links the corresponding positions of its two endpoint chains. The infinite chains distinguish the graph data from the finite ordinals and allow us to recover the original digraph after completion. Following \cite[Definition 4.1 and the proof of Theorem 4.2]{ESV}, form the disjoint union
\[
G:=\omega\sqcup(A\times\omega)\sqcup(F\times\omega).
\]
Write $a_n$ and $e_n$ for the tagged copies of $(a,n)$ and $(e,n)$, respectively, and put
\[
G_{\mathrm v}:=\{a_n:a\in A,\ n<\omega\},\qquad
G_{\mathrm e}:=\{e_n:e\in F,\ n<\omega\}.
\]
Thus $G=\omega\sqcup G_{\mathrm v}\sqcup G_{\mathrm e}$. We call the members of these three sets \emph{finite ordinal vertices}, \emph{vertex-chain vertices}, and \emph{edge-chain vertices}, respectively. For each $a\in A$ and $e\in F$, we call $a_0$ and $e_0$ the \emph{roots} of their respective chains. For $e=\{a,b\}\in F$, prescribe the fibres by
\begin{align*}
 n_E&=\{0,\ldots,n-1\},\\
 (a_n)_E&=\{n,a_{n+1}\},\\
 (e_n)_E&=\{a_n,b_n,e_{n+1}\}.
\end{align*}
\begin{figure}[ht]
\centering\includegraphics[width=.4875\linewidth]{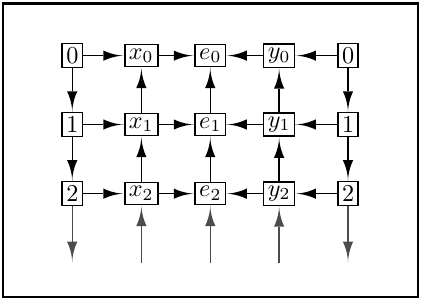}
\caption{Representation of an edge (with $x=a$ and $y=b$ in the notation of the text). Reproduced from \cite[Figure 1]{ESV}.}
\end{figure}
The ordinal fibres contain only ordinal vertices, the vertex-chain fibres contain exactly one ordinal vertex and one vertex-chain vertex, and the edge-chain fibres contain two distinct vertex-chain vertices and one edge-chain vertex. Within each type the displayed fibres determine their vertices uniquely. Thus $\mathcal G=(G,E)$ is extensional. We identify the displayed vertices with their tagged copies in the completion. We call $\mathcal G$ the \emph{seed digraph}. It is acyclic, but it has infinite descending chains. Its cardinality is less than $\kappa$. Put $\mathcal N=\mathcal G_{<\kappa}$; by Lemma~\ref{axioms}, $\mathcal N\models\mathsf{ZFC}\setminus\{\mathsf{Found}\}$.

Lemma~\ref{seed} adapts the recovery argument in the proof of \cite[Theorem 4.2]{ESV} to the transfinite completion described in Subsection~\ref{subsec:completion}.
\begin{lemma}\label{seed}
The seed $G$ is parameter-free definable in $\mathcal N$, and $\Aut(\mathcal G)\cong\Aut(\mathcal A)$.
\end{lemma}
\begin{proof}
By the fibre equations and Extensionality, external induction on $n$ identifies the ordinal vertex $n$ with the internal numeral $n^{\mathcal N}$. By Lemma~\ref{bounded}, the countable set of ordinal vertices has a representative $w$ with $w_E=\{n:n<\omega\}$. Since $0Ew$ and the internal successor of $n$ is $n+1$, $w$ is inductive. Conversely, external induction shows that every internally inductive set contains every ordinal vertex. Hence $w=\omega^{\mathcal N}$, and
\[
(\omega^{\mathcal N})_E=\{n:n<\omega\}.
\]

Using Definition~\ref{def:completion}, define $r(x)=\min\{\alpha:x\in G_\alpha\}$. The successor-stage clause gives
\[
yEx\land r(x)>0\Longrightarrow r(y)<r(x).
\]
If an externally infinite descending $E$-chain never entered $G$, its successive stages would form an infinite strictly decreasing sequence of ambient ordinals, which is impossible. Hence every such chain eventually enters $G$. Once the chain enters $G$ it remains there because the extension adds no new members to old sets.

Unlike the finite set theory setting of \cite[Definition 4.1]{ESV}, we can quantify over actual infinite sequences: by Lemma~\ref{bounded}, the internal ordered-pair graph of every external countable sequence in $N$ is represented in $\mathcal N$.

Let $\varphi_1(x)$ assert that there is a function $s$ with domain $\omega$ satisfying
\[
s(0)=x,\qquad \forall n\in\omega\ (s(n)=\{n,s(n+1)\}).
\]
We claim that
\begin{equation}\label{eq:vertex-recovery}
\mathcal N\models\varphi_1(x)\quad\Longleftrightarrow\quad
x\in\{a_0:a\in A\}.
\end{equation}
The reverse implication follows by representing $(a_n)_{n<\omega}$. For the forward implication, let $(x_n)_{n<\omega}$ be the external sequence obtained from a witness. Then
\[
x_0=x,\qquad (x_n)_E=\{n,x_{n+1}\}\quad(n<\omega).
\]
Choose $m$ with $x_m\in G$. No $x_n$ is a finite ordinal vertex, since the tail $(x_{n+i})_{i<\omega}$ would then descend through finite ordinals. Moreover,
\[
|(x_m)_E|\leq 2<|(e_k)_E|=3
\qquad(e\in F,\ k<\omega).
\]
Thus $x_m=a_k$ for some $a\in A$ and $k<\omega$, and
\[
\{m,x_{m+1}\}=\{k,a_{k+1}\}.
\]
Since $x_{m+1},a_{k+1}\notin\omega$, $m=k$ and $x_{m+1}=a_{m+1}$. Extensionality gives the backward step
\[
x_{i+1}=a_{i+1}\ \Longrightarrow\
(x_i)_E=\{i,a_{i+1}\}=(a_i)_E
\ \Longrightarrow\ x_i=a_i.
\]
For the forward step, since $x_{i+1},a_{i+1}\notin\omega$,
\[
x_i=a_i\ \Longrightarrow\ \{i,x_{i+1}\}=\{i,a_{i+1}\}
\ \Longrightarrow\ x_{i+1}=a_{i+1}.
\] Hence $x_i=a_i$ for every $i<\omega$, proving \eqref{eq:vertex-recovery} and uniqueness of the witness sequence.

Let $\theta(x,y,z)$ be the conjunction of
\[
x\ne y\land\varphi_1(x)\land\varphi_1(y)
\]
and the assertion that there exist functions $s,t,q$ with domain $\omega$ satisfying
\begin{gather*}
s(0)=x,\quad t(0)=y,\quad q(0)=z,\\
\forall n\in\omega\ \bigl(s(n)=\{n,s(n+1)\}\land
 t(n)=\{n,t(n+1)\}\bigr),\\
\forall n\in\omega\ (q(n)=\{s(n),t(n),q(n+1)\}).
\end{gather*}
We claim that
\begin{equation}\label{eq:edge-recovery}
\begin{split}
\mathcal N\models\theta(x,y,z)\quad\Longleftrightarrow\quad
&\exists a,b\in A\ \exists e\in F\ \bigl(e=\{a,b\}\\
&\land x=a_0\land y=b_0\land z=e_0\bigr).
\end{split}
\end{equation}
The reverse implication follows by representing the three prescribed sequences. For the forward implication, \eqref{eq:vertex-recovery} and uniqueness of its witnesses give distinct $a,b\in A$ and an external sequence $(z_n)_{n<\omega}$ satisfying
\[
x=a_0,\quad y=b_0,\quad z=z_0,\qquad
(z_n)_E=\{a_n,b_n,z_{n+1}\}.
\]
Choose $m$ with $z_m\in G$. The seed fibre equations give, for $u\in G$,
\[
|u_E\cap G_{\mathrm v}|=\begin{cases}
0,&u\in\omega,\\
1,&u\in G_{\mathrm v},\\
2,&u\in G_{\mathrm e}.
\end{cases}
\]
Since $\{a_m,b_m\}\subseteq(z_m)_E\cap G_{\mathrm v}$ and $a_m\ne b_m$, we have $z_m=f_k$ for some $f=\{c,d\}\in F$. In particular, $|(z_m)_E|=3$, so $z_{m+1}\notin\{a_m,b_m\}$. Comparing fibres, we have
\[
\{a_m,b_m,z_{m+1}\}=\{c_k,d_k,f_{k+1}\},
\]
Intersecting with $G_{\mathrm v}$ gives $\{a_m,b_m\}=\{c_k,d_k\}$, hence $k=m$, $f=\{a,b\}$ and $z_{m+1}=f_{m+1}$. Extensionality gives
\[
z_{i+1}=f_{i+1}\ \Longrightarrow\
(z_i)_E=\{a_i,b_i,f_{i+1}\}=(f_i)_E
\ \Longrightarrow\ z_i=f_i.
\]
Backward induction yields $z=f_0$, proving \eqref{eq:edge-recovery}. In particular, $\varphi_2(x,y):=\exists z\,\theta(x,y,z)$ recovers adjacency, and
\[
\varphi_3(z):=\exists x\exists y\,\theta(x,y,z),\qquad
\varphi_3^{\mathcal N}=\{e_0:e\in F\}.
\]

The seed consists precisely of the vertices belonging to the transitive closures of the vertex-chain and edge-chain roots, together with all the natural numbers. Consequently, $G$ is defined in $\mathcal N$ by the formula
\[
\begin{aligned}
\varphi_4(u):={}&\exists v\bigl((\varphi_1(v)\lor\varphi_3(v))
\land u\in\operatorname{TC}(\{v\})\bigr)\\
&{}\lor u\in\omega.
\end{aligned}
\]
Thus
\[
G=\{u\in N:\mathcal N\models\varphi_4(u)\}.
\]
Indeed, the seed is closed under taking members, and the transitive closure of a singleton is obtained by iterating Union through $\omega$ (without Foundation). Since $\omega^{\mathcal N}$ has exactly the standard numerals as members, these closures contain only seed vertices and include every position of the corresponding chains.

Finally, every automorphism of $\mathcal G$ fixes the finite ordinal vertices (induct on $n$). The fibre equations give
\[
G_{\mathrm v}=\{u\in G\setminus\omega:u_E\cap\omega\ne\varnothing\},
\qquad (a_n)_E\cap\omega=\{n\}.
\]
Hence every automorphism preserves $G_{\mathrm v}$ and the index of each of its members. Thus for each $n<\omega$ an automorphism $h$ induces a permutation $\pi_n$ of $A$ by $h(a_n)=(\pi_n(a))_n$. Preservation of fibres gives $\pi_{n+1}(a)=\pi_n(a)$ for every $a\in A$. Write $\pi=\pi_0$. The automorphism also preserves $G_{\mathrm e}=G\setminus(\omega\cup G_{\mathrm v})$. If $e=\{a,b\}\in F$, the fibre of $h(e_n)$ contains $(\pi(a))_n$ and $(\pi(b))_n$, so $\{\pi(a),\pi(b)\}\in F$ and $h(e_n)=(\{\pi(a),\pi(b)\})_n$. Applying the same argument to $h^{-1}$ shows that $\pi\in\Aut(\mathcal A)$. Conversely, every automorphism of $\mathcal A$ acts on these chains by its action on vertices and edges, fixing the ordinal vertices. These operations are inverse.
\end{proof}

\begin{remark}
For the automorphism-group conclusion, setwise invariance of the seed suffices. The recovery argument of Lemma~\ref{seed} can be carried out externally, using external infinite sequences and transitive closures in place of their internal counterparts. In this form it establishes invariance at every ordinal stage of the Tzouvaras completion of the coding digraph, without any cardinal hypothesis.
\end{remark}

\begin{proof}[Proof of Theorem~\ref{inacc}]
By Proposition~\ref{completion-aut}, every automorphism of the seed extends uniquely to the completion (the corresponding finitary universal property appears in \cite[Remark 3.14(a)]{ESV}). By Lemma~\ref{seed}, every automorphism of the completion preserves the seed. Since $\mathcal G_{<\kappa}=\mathcal G_\kappa$, restriction and extension are inverse isomorphisms of abstract groups, and
\[
\Aut(\mathcal N)\cong\Aut(\mathcal G)\cong\Aut(\mathcal A)\cong H.
\]
The cardinality and axioms other than Foundation follow from Lemmas~\ref{bounded} and \ref{axioms}. For any $a\in A$, the set representing $\{a_n:n<\omega\}$ witnesses the failure of Foundation. Proposition~\ref{wf-completion} gives the assertion about the well-founded part. Moreover, the completion is acyclic: a cycle meeting a new vertex would contradict descent of construction ranks, and the seed itself is acyclic.
\end{proof}

\begin{remark}
The proper-class construction of Appendix~\ref{sec:class-completions}, applied to the coding digraph of Lemma~\ref{seed}, also yields a proper-class model for each set-sized group. In Kelley--Morse class theory, external recovery and transfinite recursion give mutually inverse restriction and extension correspondences between class automorphisms of the completion and automorphisms of the seed. No collection of proper classes is required as an object of the class theory. Theorem~\ref{inacc}, in contrast, concerns set structures.

A strongly inaccessible cardinal is a sufficient hypothesis for the set-sized construction; it remains to determine whether one can prove that every set-sized group is isomorphic to the automorphism group of a set-sized model of $\mathsf{ZFC}^{-\mathsf{Found}}$ in \textup{\textsf{ZFC}} alone. If there is a proper class of strongly inaccessible cardinals, Theorem~\ref{inacc} applies to every group. 
\end{remark}

\section{Questions}\label{sec:open-question}

Theorems A and B illustrate the contrasting roles of Choice
and Foundation in determining the automorphism groups of
models of set theory. By Theorem A(1),
\[
\Grpr\subsetneq\G(\ZF)\subseteq\Grp.
\]
Theorem B shows that, under a suitable cardinal hypothesis,
every group can be realized in the absence of Foundation,
even in the presence of Choice.
This suggests the following question concerning the
effect of dispensing with Choice while retaining Foundation.

\textbf{Question 1.} Is $\G(\ZF)=\Grp$?

By Theorem A(1),
\[
\G(\mathsf{ZFC})=\G(\ZF+\mathsf{LM})
=\G(\ZF+\mathsf{V=OD})=\Grpr.
\]
Moreover, Corollary~\ref{od-spectrum} establishes that
this equality persists for every consistent extension
of $\ZF+\mathsf{V=OD}$ in $\Lset$,
whereas Proposition~\ref{faithful} and Theorem~\ref{orderability}
ensure that $\G(T)\subseteq\Grpr$ for every consistent
extension $T$ of $\mathsf{ZFC}$ or $\ZF+\mathsf{LM}$.

\textbf{Question 2.} Is $\G(T)=\Grpr$ for every consistent
extension $T$ of $\mathsf{ZFC}$ or of $\ZF+\mathsf{LM}$
in $\Lset$? In particular, is
\[
\G(\ZF+\mathsf{LM}+\neg\mathsf{AC})=\Grpr?
\]

In contrast to the realization problems posed in
Questions 1 and 2, the Ehrenfeucht--Mostowski theorem
ensures that, for every first-order theory $T$ with
an infinite model, every right-orderable group embeds
into the automorphism group of some model of $T$
(see Section~\ref{sec:intro}).

\appendix
\section{Proper-class completions and descending membership chains}\label{sec:class-completions}

The completion construction is of interest independently of the study of automorphisms. Tzouvaras \cite[Section~4]{Tz} carries out a proper-class construction using Global Choice. Our variant, based on tagging and local set-sized recursion, can be carried out over $\ZF$ alone. We then use it to address the question of external descending membership chains in ordinal-standard models raised in \cite{MO}.

\begin{definition}
Let $\mathcal M\models\ZF$ and let $\mathcal G=(G,E)$ be a parametrically $\mathcal M$-definable class digraph. For $s\in M$, write $s^{\mathcal M}=\{x\in M:x\in^{\mathcal M}s\}$. \begin{enumerate}[label=(\alph*)]
\item $\mathcal G$ is \emph{set-like in $\mathcal M$} if for every $a\in G$ there is an $s\in M$ with $a_E=s^{\mathcal M}$.
\item $\mathcal G$ is \emph{full in $\mathcal M$} if for every $s\in M$ with $s^{\mathcal M}\subseteq G$ there is an $a\in G$ with $a_E=s^{\mathcal M}$.
\item The \emph{well-founded part of $\mathcal G$ in
the sense of $\mathcal M$} consists of those vertices
whose downward closures, as computed in
$\mathcal M$, are well-founded in $\mathcal M$.
\end{enumerate}
\end{definition}

\begin{theorem}\label{class-completion}
Suppose $\mathcal M\models\ZF$ and $\mathcal G$ is a parametrically $\mathcal M$-definable extensional class digraph that is set-like in $\mathcal M$. Then $\mathcal G$ admits a parametrically $\mathcal M$-definable end embedding into a parametrically $\mathcal M$-definable extensional digraph $\mathcal N=(N,F)$ such that:
\begin{enumerate}[label=(\alph*)]
\item $\mathcal N$ is set-like and full in $\mathcal M$.
\item $\mathcal N\models\ZF\setminus\{\mathsf{Found}\}$, and if $\mathcal M\models\mathsf{AC}$, then $\mathcal N\models\mathsf{AC}$.
\item There is a canonical parametrically $\mathcal M$-definable end embedding $j:\mathcal M\longrightarrow\mathcal N$ whose image is precisely the well-founded part of $\mathcal N$ in the sense of $\mathcal M$.
\end{enumerate}
\end{theorem}
\begin{proof}
Work within $\mathcal M$. Let $(G_0,E_0)$ be the tagged copy of the seed and let $D_0$ be the definable class of sets represented by its fibres. Define $N$ by the recursive condition
\[
u\in N\iff u\in G_0\lor\exists s\,(u=\langle1,s\rangle\land s\subseteq N\land s\notin D_0).
\]
Every member of $s$ has smaller rank than $\langle1,s\rangle$. Thus the condition can be implemented by set-sized well-founded recursion on the transitive closure of $\{u\}$. Uniqueness of the local recursions makes their answers coherent, giving a single definition of $N$ with the parameters of the seed. Retain the old fibres and give each new vertex $\langle1,s\rangle$ the fibre $s$. Extensionality follows from tagging and the exclusion of $D_0$. Every set contained in $N$ is represented, either by an old vertex or by its tag. This establishes (a), with the tagging map giving the required end embedding of $\mathcal G$.

For (b), the verification of the axioms is as in Lemma~\ref{axioms}, without cardinal bounds. For Separation and Replacement translate each individual formula into the language of $\mathcal M$ using the definitions of $N$ and $F$. For Collection, Collection in $\mathcal M$ applied to the translated formula supplies a set of witnesses, whose representative witnesses the corresponding instance in $\mathcal N$. Infinity follows by the same recursive construction of the finite ordinals. Choice in $\mathcal M$, when available, supplies well-orderings of the fibres. This is the full representation property of Rieger's theorem (see \cite[Appendix B]{Ac}).

For (c), write $[s]$ for the unique representative of a set $s\subseteq N$. Ordinary set recursion defines
\[
j(a)=[\{j(b):b\in a\}].
\]
Induction on rank shows that this is an isomorphism onto a downward-closed submodel of $\mathcal N$. Conversely, every vertex with well-founded downward closure belongs to its range by the Mostowski collapse. The downward closure is a set, obtained by iterating the operation of taking $F$-predecessors through $\omega$, using Replacement and Union. Every ordinal of $\mathcal N$ belongs to this part: the ordering on its fibre is a well-order in $\mathcal M$, since every subset in $\mathcal M$ is represented. Conversely, the images of the ordinals of $\mathcal M$ are ordinals of $\mathcal N$.
\end{proof}

Note that $\mathcal N$ satisfies Foundation if and only
if the seed digraph $\mathcal G$ is well-founded in the
sense of $\mathcal M$, in which case the canonical
end embedding $j$ is an isomorphism.

In particular, if $\mathcal M$ is standard, then $\mathcal N$ is ordinal-standard of the same ordinal height, and $j$ identifies $\mathcal M$ with the externally well-founded part of $\mathcal N$.

\begin{remark}
The recursive definition of the completion itself works even when the seed is not set-like in $\mathcal M$, by taking $D_0$ to consist of those sets of $\mathcal M$ that are represented. Set-likeness is needed for the general assertion about the axioms. For example, working within $\mathcal M$, take a reverse ordinal chain indexed by all ordinals. Removing one intermediate vertex from a fibre gives a proper-class subclass which is neither an original fibre nor a set; the corresponding instance of Separation fails in the completion. This is a limitation of the construction described in \cite{MO}. Nevertheless, as shown in Proposition~\ref{chain-existence},
a variant of that construction yields a positive answer
to the question raised in \cite{MO}.
\end{remark}

Tzouvaras \cite[Theorem~3.5]{Tz} rules out a descending path of length the entire class of ordinals under his class-theoretic hypotheses. The question in \cite{MO}, however, concerns an externally given chain, which need not be available internally. The following two propositions address that question.

\begin{definition}\label{def:descending-chains-paths}
Let $\mathcal N=(N,F)$ be an $\Lset$-structure and let $\beta$ be an ordinal.
\begin{enumerate}[label=(\alph*)]
\item An \emph{external descending membership chain} of length $\beta$ is an injective sequence $(a_\xi:\xi<\beta)$ of elements of $N$ such that $a_\eta F a_\xi$ whenever $\xi<\eta<\beta$.
\item Following Tzouvaras \cite[Definition~3.1]{Tz}, an \emph{external descending membership path} of length $\beta$ is an injective sequence $(a_\xi:\xi<\beta)$ of elements of $N$ such that $a_{\xi+1} F a_\xi$ whenever $\xi+1<\beta$, and for every nonzero limit ordinal $\lambda<\beta$ there is a $\gamma<\lambda$ such that
\[
a_\lambda F a_\xi
\quad\text{whenever}\quad \gamma\leq\xi<\lambda.
\]
\end{enumerate}
\end{definition}

Every external descending membership chain is a path, but not conversely.

\begin{proposition}\label{chain-obstruction}
If an ordinal-standard model of $\mathsf{ZFC}^{-\mathsf{Found}}$ of ordinal height $\beta$ admits an external descending membership chain of length $\beta$ in the above sense, then $\beta$ is not a cardinal (i.e., an initial ordinal).
\end{proposition}
\begin{proof}
By Choice every internal set has an internal bijection with an ordinal, which externally gives a bijection between its fibre and some ordinal less than $\beta$. The fibre of the first vertex of the chain contains all the remaining vertices. Since $\beta$ is an infinite limit ordinal, there is therefore an $\alpha<\beta$ with $|\alpha|=|\beta|$. Hence $\beta$ is not a cardinal.
\end{proof}
\begin{remark}
If the ordinal height $\beta$ in
Proposition~\ref{chain-obstruction} is an uncountable
regular cardinal, then even an external descending
membership path of length $\beta$ cannot exist.
Indeed, by Fodor's lemma, the regressive function $f$ obtained by choosing a witness $f(\lambda)<\lambda$ to the path condition at each nonzero limit ordinal $\lambda<\beta$ is constant, say with value
$\gamma$, on a stationary subset of $\beta$.
Consequently, the fibre of $a_\gamma$ has cardinality
at least $\beta$, contradicting Choice and
ordinal-standardness as in the proof of
Proposition~\ref{chain-obstruction}.
This is an external counterpart of
Tzouvaras's Theorem~3.5.
\end{remark}

\begin{proposition}\label{chain-existence}
Suppose $\mathcal M$ is a standard model of $\mathsf{ZFC}$ of countable ordinal height $\beta$. Then there is an acyclic, parametrically $\mathcal M$-definable model $\mathcal N$ of $\mathsf{ZFC}^{-\mathsf{Found}}$ such that:
\begin{enumerate}[label=(\alph*)]
\item $\mathcal M$ is canonically isomorphic to the externally well-founded part of $\mathcal N$.
\item $\mathcal N$ admits an external descending membership chain of length $\beta$.
\end{enumerate}
\end{proposition}
\begin{proof}
Within $\mathcal M$ apply Theorem~\ref{class-completion} to the reverse rational order, with $qEr$ iff $r<q$. Its fibres are distinct sets, and it is acyclic. The set representing the range of a strictly increasing rational sequence witnesses the failure of Foundation. The resulting completion remains acyclic by descent of ranks in $\mathcal M$ outside the seed. Part (c) of Theorem~\ref{class-completion} gives (a). Externally every countable linear order embeds in the rational order, which is standard since $\mathcal M$ is standard; an embedding of $\beta$ therefore gives the desired descending chain. It need not be present in $\mathcal M$.
\end{proof}

\begin{remark}
By a theorem of Friedman \cite[Theorem~3.4]{Friedman1973}, if $\mathcal M$ is a countable $\omega$-nonstandard
model  of $\ZF$, then $(\operatorname{Ord}^{\mathcal M},\in^{\mathcal M})$ is isomorphic to the sum
of an initial segment and a $\mathbb Q$-indexed
lexicographic sum of nonempty countable linear orders
determined by the standard system of $\mathcal M$,
where $\mathbb Q$ denotes the order type of the rationals.
In particular, there is a copy of $\mathbb Q$ in
$(\operatorname{Ord}^{\mathcal M},\in^{\mathcal M})$.
Since every countable linear order embeds into
$\mathbb Q$, there is an external order embedding
\[
(\operatorname{Ord}^{\mathcal M},\ni^{\mathcal M})
\longrightarrow
(\operatorname{Ord}^{\mathcal M},\in^{\mathcal M}),
\]
where $\ni^{\mathcal M}$ denotes the reversal of
$\in^{\mathcal M}$. 
\end{remark}

By carrying out Tzouvaras's original construction \cite[Theorem~4.2(a)]{Tz} over $V_\kappa$, where $\kappa$ is strongly inaccessible, with an external global well-order of $V_\kappa$ of order type $\kappa$, we obtain an acyclic ordinal-standard model of $\mathsf{ZFC}^{-\mathsf{Found}}$ of ordinal height $\kappa$ that contains an internally represented strictly descending membership chain of every internal ordinal length. By Proposition~\ref{chain-obstruction}, however, no such model admits an external descending membership chain of length $\kappa$. This illustrates the distinction between the existence of arbitrarily long internal descending membership chains and the existence of an external descending membership chain indexed by the entire ordinal height of the model.

\section*{Acknowledgment}By 2012, I had developed the central ideas and proof strategies of this paper, though they remained unpublished. Of course, this involved standing on the shoulders of many giants, as reflected in the way various disparate ideas are combined here so as to see a little farther.

ChatGPT Edu (courtesy of the University of Gothenburg) came to my aid as an indefatigable and knowledgeable scribe to assist with drafting the exposition, especially with the routine---but sometimes laborious---details, all of which I have checked, corrected, and edited (in more rounds than I would have imagined). 


\begin{thebibliography}{99}
\raggedright
\bibitem{AH} F. G. Abramson and L. A. Harrington, Models without indiscernibles, \emph{Journal of Symbolic Logic} \textbf{43} (1978), no.~3, 572--600.
\bibitem{Ack} W. Ackermann, Zur Widerspruchsfreiheit der Zahlentheorie, \emph{Mathematische Annalen} \textbf{117} (1940), 162--194.
\bibitem{Ac} P. Aczel, \emph{Non-Well-Founded Sets}, CSLI Lecture Notes 14, Stanford, 1988.
\bibitem{An1} D. A. Anapolitanos, Automorphisms of finite order, \emph{Zeitschrift f\"ur Mathematische Logik und Grundlagen der Mathematik} \textbf{25} (1979), no.~6, 565--575 (MR0549475).
\bibitem{An2} D. A. Anapolitanos, Automorphisms of finite order and the axiom of choice, \emph{Bulletin de la Soci\'et\'e Math\'ematique de Belgique, S\'erie B} \textbf{31} (1979), no.~1, 39--44 (MR0592658).
\bibitem{Babai} L. Babai, Infinite digraphs with given regular automorphism groups, \emph{Journal of Combinatorial Theory, Series B} \textbf{25} (1978), 26--46.
\bibitem{Cohen} P. J. Cohen, Automorphisms of set theory, in \emph{Proceedings of the Tarski Symposium}, Proceedings of Symposia in Pure Mathematics 25, American Mathematical Society, 1974, pp. 325--330.\bibitem{Cohn}
P.~M. Cohn,
Groups of order automorphisms of ordered sets,
\emph{Mathematika} \textbf{4} (1957), 41--50.

\bibitem{Conrad}
P. Conrad,
Right-ordered groups,
\emph{Michigan Mathematical Journal}
\textbf{6} (1959), 267--275.
\bibitem{deGroot} J. de Groot, Groups represented by homeomorphism groups I, \emph{Mathematische Annalen} \textbf{138} (1959), 80--102.
\bibitem{Bogota} A. Enayat, Analogues of the MacDowell--Specker theorem for set theory, in X. Caicedo and C. H. Montenegro (eds.), \emph{Models, Algebras and Proofs}, Marcel Dekker, 1998, pp. 25--50; revised version, September 14, 2024.
\bibitem{EM} A. Enayat, Automorphisms, Mahlo cardinals, and NFU, in A. Enayat and R. Kossak (eds.), \emph{Nonstandard Models of Arithmetic and Set Theory}, Contemporary Mathematics 361, American Mathematical Society, 2004, pp. 37--59.
\bibitem{LM} A. Enayat, On the Leibniz--Mycielski axiom in set theory, \emph{Fundamenta Mathematicae} \textbf{181} (2004), 215--231.
\bibitem{Paris} A. Enayat, Models of set theory with definable ordinals, \emph{Archive for Mathematical Logic} \textbf{44} (2005), 363--385.
\bibitem{Karagila} A.~Karagila, \emph{The Morris model}, Proc. Amer. Math. Soc. \textbf{148} (2020), 1311--1323.
\bibitem{Morris} D. B. Morris, Existence of models of ZF with involutions, preliminary report, \emph{Notices of the American Mathematical Society} \textbf{16} (1969), no.~6, 980, abstract 69T-E72.
\bibitem{MorrisThesis} D.~B. Morris, \emph{Adding total indiscernibles to models of set theory}, Ph.D. thesis, University of Wisconsin--Madison, 1970.
\bibitem{MO} A. Enayat, Answer to ``Ill-founded models of set theory with well-founded ordinals'', MathOverflow, 1 January 2016, \url{https://mathoverflow.net/a/227423}.
\bibitem{Dead} A. Enayat, Models of set theory: Extensions and dead-ends, \emph{Journal of Symbolic Logic} (2026), doi:10.1017/jsl.2026.10199.
\bibitem{EKM} A. Enayat, M. Kaufmann and Z. McKenzie, Largest initial segments pointwise fixed by automorphisms of models of set theory, \emph{Archive for Mathematical Logic} (2017), doi:10.1007/s00153-017-0582-3.
\bibitem{EL} A. Enayat and M. \L e\l yk, Categoricity-like properties in the first order realm, \emph{Journal for the Philosophy of Mathematics} \textbf{1} (2024), 63--98.
\bibitem{ESV} A. Enayat, J. H. Schmerl and A. Visser, $\omega$-models of finite set theory, in \emph{Set Theory, Arithmetic, and Foundations of Mathematics: Theorems, Philosophies}, Lecture Notes in Logic, Cambridge University Press, 2011.
\bibitem{Forster} T. E. Forster, \emph{Set Theory with a Universal Set: Exploring an Untyped Universe}, 2nd ed., Oxford University Press, 1995.
\bibitem{Friedman1973}
H. Friedman, Countable models of set theories,
in A. R. D. Mathias and H. Rogers (eds.),
\emph{Cambridge Summer School in Mathematical Logic},
Lecture Notes in Mathematics 337, Springer-Verlag,
1973, pp. 539--573.

\bibitem{Frucht} R. Frucht, Herstellung von Graphen mit vorgegebener abstrakter Gruppe, \emph{Compositio Mathematica} \textbf{6} (1939), 239--250.
\bibitem{Gaifman1974} H. Gaifman, Operations on relational structures, functors and classes I, in \emph{Proceedings of the Tarski Symposium}, Proceedings of Symposia in Pure Mathematics 25, American Mathematical Society, 1974, pp. 21--39.
\bibitem{Gaifman} H. Gaifman, Models and types of Peano's arithmetic, \emph{Annals of Mathematical Logic} \textbf{9} (1976), 223--306.
\bibitem{Hodges} W. Hodges, \emph{Model Theory}, Encyclopedia of Mathematics and its Applications 42, Cambridge University Press, 1993.
\bibitem{HodgesReview} W. Hodges, Review of D. A. Anapolitanos, Automorphisms of finite order and the axiom of choice, \emph{Mathematical Reviews}, MR0592658.
\bibitem{Kaufmann} M. Kaufmann, Blunt and topless end extensions of models of set theory, \emph{Journal of Symbolic Logic} \textbf{48} (1983), 1053--1073.
\bibitem{KW} R. Kaye and T. Wong, On interpretations of arithmetic and set theory, \emph{Notre Dame Journal of Formal Logic} \textbf{48} (2007), 497--510.
\bibitem{KS} R. Kossak and J. H. Schmerl, \emph{The Structure of Models of Peano Arithmetic}, Oxford Logic Guides 50, Oxford University Press, 2006.
\bibitem{Montague} R. Montague, Semantical closure and non-finite axiomatizability I, in \emph{Infinitistic Methods}, Pergamon Press, 1961, pp. 45--69.
\bibitem{PillayShelah} A. Pillay and S. Shelah, Classification theory over a predicate. I, \emph{Notre Dame Journal of Formal Logic} \textbf{26} (1985), 361--376.
\bibitem{Pillay2026} A. Pillay, \emph{Remarks on relative categoricity}, arXiv:2602.05866v1, 2026.
\bibitem{Sabidussi} G. Sabidussi, Graphs with given infinite group, \emph{Monatshefte f\"ur Mathematik} \textbf{64} (1960), 64--67.
\bibitem{Schmerl} J. H. Schmerl, Automorphism groups of models of Peano arithmetic, \emph{Journal of Symbolic Logic} \textbf{67} (2002), no.~4, 1249--1264.
\bibitem{Togha} A. A. Togha, Automorphisms of countable recursively saturated models of set theory, in A. Enayat and R. Kossak (eds.), \emph{Nonstandard Models of Arithmetic and Set Theory}, Contemporary Mathematics 361, American Mathematical Society, 2004, pp. 163--167.
\bibitem{Tz} A. Tzouvaras, Non-circular non-well-founded set universes, \emph{Mathematical Logic Quarterly} \textbf{39} (1993), 454--460.
\bibitem{Usvyatsov} A. Usvyatsov, \emph{On the existence property over a predicate}, arXiv:2502.20236, 2025.
\bibitem{Visser2006} A. Visser, Categories of theories and interpretations, in \emph{Logic in Tehran}, Lecture Notes in Logic 26, Association for Symbolic Logic, 2006, pp. 284--341.
\end{thebibliography}
\end{document}